\documentclass[a4paper,10pt,leqno]{article}

\usepackage[T1]{fontenc}
\usepackage{lmodern}
\usepackage{microtype}

\usepackage{amsmath,amssymb,amsfonts,amsthm}
\usepackage{mathtools}
\usepackage{mathrsfs}
\usepackage{bm}
\usepackage{dsfont}

\usepackage{graphicx}
\usepackage{float}
\usepackage{booktabs,multirow}

\usepackage{enumitem}

\usepackage[margin=1.2in]{geometry}

\usepackage{xcolor}
\usepackage[colorlinks=true,linkcolor=black,citecolor=black,urlcolor=blue]{hyperref}

\usepackage{tikz}
\usetikzlibrary{arrows.meta}

\allowdisplaybreaks
\mathtoolsset{showonlyrefs}

\numberwithin{equation}{section}

\theoremstyle{plain}
\newtheorem{thm}{Theorem}[section]

\newtheorem{lem}[thm]{Lemma}
\newtheorem{prop}[thm]{Proposition}
\newtheorem{dfn}[thm]{Definition}

\theoremstyle{remark}
\newtheorem{rmk}[thm]{Remark}

\newtheoremstyle{main}
  {6pt}{6pt}
  {\itshape}{}
  {\scshape}{.}{0.6em}{}
\theoremstyle{main}
\newtheorem{mainthm}{Theorem}
\newtheorem{maincor}{Corollary}

\newcommand{\e}{\varepsilon}

\DeclareMathOperator{\Fix}{Fix}
\DeclareMathOperator{\diag}{diag}
\DeclareMathOperator{\Graph}{Graph}
\DeclareMathOperator{\grad}{grad}
\DeclareMathOperator{\Hess}{Hess}

\newcommand{\LC}{\mathcal{LC}}

\newcommand{\U}{\mathbb{U}}
\newcommand{\R}{\mathbb{R}}
\newcommand{\C}{\mathbb{C}}
\newcommand{\Z}{\mathbb{Z}}
\newcommand{\N}{\mathbb{N}}
\newcommand{\Id}{\mathrm{Id}}
\newcommand{\Sp}{\mathrm{Sp}}

\newcommand{\iMorT}[1]{\iota^{T}\!\left(#1\right)}
\newcommand{\iMorE}[1]{\iota^{E}\!\left(#1\right)}
\newcommand{\iMorD}[1]{\iota_D\!\left(#1\right)}
\newcommand{\iCLM}{\iota^{\scriptscriptstyle{\mathrm{CLM}}}}

\title{
Non-minimality and instability of brake orbits 
for natural Lagrangians on Riemannian manifolds
}

\author{
Luca Asselle\thanks{L.A. is partially supported by the DFG-Grant 566804407 “Symplectic Dynamics, Celestial Mechanics and magnetism”} \and
Xijun Hu\thanks{X.H. is partially  supported by the National Natural Science Foundation of China (No. 12521001) and Taishan Scholars Climbing Program of Shandong (TSPD20240802)}
\and
Alessandro Portaluri\thanks{A.P. is partially supported by Tamkeen Portaluri 2025--26 Faculty Research Funds, Tamkeen NYU Abu Dhabi, and by GNAMPA--INDAM, Italy.}
\and
Li Wu\thanks{L.W. is partially supported by the National Natural Science Foundation of China (NSFC No.\ 12171281).}
}

\date{\today}

\begin{document}
\maketitle

\begin{abstract}
We investigate minimality and stability of periodic brake orbits in natural Lagrangian systems on smooth Riemannian manifolds. We prove that every non-constant periodic brake orbit is not a minimizer of the fixed-time action, for any conormal boundary condition. Under an orbit-cylinder hypothesis, its Morse index strictly increases in the free-time setting.

As a consequence, strongly nondegenerate brake orbits fail to be linearly stable under a dimensional condition; in dimension at least three, nondegenerate mountain-pass brake orbits are spectrally unstable when the monodromy is semisimple.

The key ingredient is a local index contribution at each brake instant. Using Seifert collar coordinates near the Hill boundary, we reduce the normal dynamics to a one-dimensional model, exhibiting a degeneracy inherent to brake symmetry.

We illustrate the results by explicit Morse index computations for the planar anisotropic oscillator, the planar pendulum, and the planar Kepler problem; in the Kepler case, the ejection--collision orbit is treated via cotangent-lift Levi--Civita--Lissajous regularization.


\end{abstract}


\section{Introduction, description of the problem and main results}

\emph{Brake orbits}, introduced by Seifert~\cite{Sei48}, are periodic solutions of conservative second-order systems that are \emph{time-reversal symmetric}: $q(-t)=q(t)$ and $\dot q(-t)=-\dot q(t)$. The instants where $\dot q(t)=0$ are \emph{brake points}: the motion stops, reverses, and retraces the same geometric arc. They arise naturally in time-reversible models and are central in Celestial Mechanics.

Their existence, multiplicity, and stability have been widely studied via variational and symplectic methods. Multiplicity results based on Maslov-type index theory were obtained by Liu and Zhang~\cite{LZ14}; instability of symmetric periodic solutions was analyzed by Ure\~na~\cite{Ure06,Ure13,Ure18}; and Hu, Wu, and Yang~\cite{HWY20} established sharp relations between Morse index and linear stability.

We prove that, under general assumptions, a $T$-periodic brake orbit $\gamma$ of a natural Lagrangian system on a Riemannian manifold:
\begin{itemize}
\item[(O1)] is never a minimizer, neither for the fixed-time problem (under any self-adjoint boundary condition) nor for the free-time problem;
\item[(O2)] is not linearly stable whenever $\dim M - 2\,\iMorT{\gamma} \ge 1$,
\end{itemize}
where $\iMorT{\gamma}$ is the fixed-time Morse index, i.e.\ the dimension of the negative spectral subspace of the second variation of the fixed-time action at $\gamma$.


\subsection{Description of the problem and main results}

Let $(M,g)$ be a smooth Riemannian manifold and $V\in\mathscr C^{2}(M)$. Consider the natural Lagrangian
\begin{equation}\label{eq:natural-L}
L:TM\to\R,\qquad L(q,v)=\dfrac12\,g_q(v,v)-V(q).
\end{equation}
For $k\in\R$ define the Hill's region and the Hill's boundary as
\[
\mathcal H_k=\{q\in M:V(q)\le k\},\qquad 
\partial\mathcal H_k=\{q\in M:V(q)=k\}.
\]
Assume that $k$ is a regular value of $V$,
\begin{equation}\label{eq:no-flat}
dV(q)\neq0 \quad \text{for all } q\in\partial\mathcal H_k,
\end{equation}
so that $\partial\mathcal H_k$ is a $\mathscr C^{2}$ hypersurface (equivalently, $\nabla V\neq0$ on $\partial\mathcal H_k$).

By Maupertuis' principle, energy-$k$ trajectories are reparametrized geodesics of the Jacobi--Maupertuis metric
\[
g_k(q)(\dot q,\dot q)=2\big(k-V(q)\big)\,g_q(\dot q,\dot q),
\qquad q\in\mathcal H_k^\circ,
\]
which degenerates on $\partial\mathcal H_k$ and produces index contributions near brake instants (cf.\ \cite{Mon14}).

A $T$-periodic brake solution at energy $k$ is a pair $(T,q)$ such that
\begin{equation}\label{eq:brake-orbit-Lagrangian}
\begin{cases}
\dfrac{D}{dt}\dot q(t)+\nabla V(q(t))=0,
\qquad E(q(t),\dot q(t))\equiv k,\\[4pt]
q(-t)=q(t),\ \dot q(-t)=-\dot q(t),
\qquad q(t+T)=q(t),
\end{cases}
\quad \forall t\in\R,
\end{equation}
where $D/dt$ is the covariant derivative and $E(q,v)=\dfrac12|v|_g^{2}+V(q)$.
Via the Legendre transform, the Hamiltonian on $T^{*}M$ is
\[
H(q,p)=\dfrac12\,g_q^{-1}(p,p)+V(q),
\]
and time reversal is $S(q,p)=(q,-p)$, with $\Fix(S)=M\times\{0\}$; hence brake points are exactly the instants where $p=0$.

\begin{rmk}[Brake instants and Hill boundary]
In dimension one, a brake instant is a turning point. In higher dimension, $\partial\mathcal H_k$ is a hypersurface and, near a brake instant, Seifert collar coordinates show that the tangential motion is higher order, so the orbit approaches $\partial\mathcal H_k$ essentially along the normal direction. This yields the local ``throwing-ball'' behavior underlying the index estimate. We focus on $T$-periodic brake solutions touching $\partial\mathcal H_k$; then the orbit meets the Hill boundary at two brake instants per period and is geometrically a chord of $\partial\mathcal H_k$.
\end{rmk}

\subsection*{Variational setting and Morse indices}

Fixed-period solutions are critical points of the fixed-time action functional
\begin{equation}\label{eq:fixed-time-action-intro}
\mathcal A^T(q)=\int_0^T L\bigl(q(t),\dot q(t)\bigr)\,dt
\end{equation}
defined on the Hilbert manifold $H^1_{\mathrm{per}}([0,T],M)$.

\begin{dfn}[Fixed-time Morse indices]\label{def:fixed-time-Morse-intro}
Let $\gamma$ be a non-constant $T$-periodic solution of the Euler--Lagrange
equation
\begin{equation}\label{eq:EL}
\dfrac{D}{dt} \dot q(t) + \nabla V\big(q(t)\big)=0, \qquad t\in(0,T].
\end{equation}
The \emph{fixed-time Morse index} $\iMorT{\gamma}$ is the dimension of the
maximal subspace of $T_\gamma H^1_{\mathrm{per}}([0,T],M)$ on which
$d^2\mathcal A^T(\gamma)$ is negative definite. The \emph{Dirichlet fixed-time
Morse index} $\iMorD{\gamma}$ is the dimension of the maximal subspace of
$T_\gamma H^1_{(p,q)}([0,T],M)$ (with fixed endpoints $p,q\in M$) on which
$d^2\mathcal A^T(\gamma)$ is negative definite.
\end{dfn}

\subsubsection*{Orbit cylinders, (non)degeneracy, and the free-time index}

Fixed-time periodic orbits are characterized variationally as critical points
of the free-period action functional. Given $x:\R/\Z\to M$ and $T>0$, set
$q(t):=x(t/T)$; then
\[
\int_0^T L(q(t),\dot q(t))\,dt
=
T\int_0^1 L\bigl(x(s),x'(s)/T\bigr)\,ds.
\]
Fixing the energy level $k$ leads to the augmented functional
\begin{equation}\label{eq:free-period-functional-intro}
\mathcal S_k(x,T)=T\int_0^1\bigl[L(x(s),x'(s)/T)+k\bigr]\,ds,
\qquad (x,T)\in H^1(\R/\Z,M)\times(0,+\infty).
\end{equation}
Critical points $(x,T)$ of $\mathcal S_k$ correspond to $T$-periodic solutions of
\eqref{eq:EL} with energy $k$.

\begin{dfn}[Free-time Morse index]\label{def:free-time -Morse-intro}
Let $(x,T)$ be a non-constant critical point of $\mathcal S_k$. The
\emph{free-time Morse index} $\iMorE{x}$ is the dimension of the maximal
subspace of $T_{(x,T)}\!\big(H^1(\R/\Z,M)\times(0,+\infty)\big)$ on which
$d^2\mathcal S_k(x,T)$ is negative definite.
\end{dfn}

For natural $\mathscr C^2$ Lagrangians these indices are finite (cf.\
\cite{Abb13,PWY22}).

\begin{dfn}[Orbit cylinder]\label{def:orbit-cylinder}
We say that $(x,T)$ (or the corresponding periodic solution $\gamma$) admits an
\emph{orbit cylinder} if there exist $\varepsilon>0$ and a $\mathscr \mathscr C^1$ family
\[
s\longmapsto (x_{k+s},T_{k+s})\in H^1(\R/\Z,M)\times(0,+\infty),
\qquad |s|<\varepsilon,
\]
such that $(x_{k+s},T_{k+s})$ is a critical point of $\mathcal S_{k+s}$ for every
$|s|<\varepsilon$ and $(x_k,T_k)=(x,T)$.
\end{dfn}

\begin{dfn}[Non-degenerate and degenerate cylinders]\label{def:nondeg-cylinder}
An orbit cylinder is \emph{non-degenerate} at $k$ if $T'(k)\neq 0$ and
\emph{degenerate} at $k$ if $T'(k)=0$.
\end{dfn}


Assume that $\gamma$ belongs to an orbit cylinder near $k$ in the sense of
Definition~\ref{def:orbit-cylinder}. Then the fixed-time and free-time Morse
indices satisfy
\begin{equation}\label{eq:index-diff-cylinder}
\iMorE{\gamma}=\iMorT{\gamma}+C(k),
\qquad C(k)\in\{0,1\}.
\end{equation}
Moreover, the correction term is determined by the monotonicity of the period:
\begin{equation}\label{eq:Ck-sign}
C(k)=
\begin{cases}
1, & \text{if } T'(k)\ge 0,\\[0.2em]
0, & \text{if } T'(k)<0.
\end{cases}
\end{equation}
Here $C(k)$ is the contribution coming from the ``period direction'' in the
free-period Hessian; see \cite{PWY22} for the general Tonelli setting.

\begin{rmk}[Degenerate cylinders]\label{rmk:degenerate-cylinder}
If $T'(k)=0$, the cylinder is degenerate and the free-period critical point
$(x,T)$ may have positive nullity. Nevertheless,
\eqref{eq:index-diff-cylinder}--\eqref{eq:Ck-sign} still give $C(k)=1$ under the
monotonicity assumption $T'(k)\ge 0$, hence
\[
\iMorE{\gamma}=\iMorT{\gamma}+1.
\]
Thus degeneracy affects the kernel of the free-period Hessian, but not the index
inequality used in the main theorem.
\end{rmk}

\begin{mainthm}\label{thm:main-1-intro}
Let $(M,g)$ be a Riemannian manifold and let $L$ be the natural Lagrangian in
\eqref{eq:natural-L}. Assume that $k$ is a regular value of $V$ and that
$\partial\mathcal H_k\neq\emptyset$. Let $\gamma$ be a non-constant $T$-periodic
brake orbit of energy $k$. Then
\[
\iMorT{\gamma} \ge \iMorD{\gamma}\ge 1,
\]
and hence $\gamma$ cannot be a fixed-time minimizer of $\mathcal A^T$, regardless
of the imposed conormal boundary conditions.

Assume in addition that $\gamma$ lies on an orbit cylinder
$\{(\gamma_{k+s},T(k+s))\}_{|s|<\varepsilon}$ and that $T'(k)\ge 0$. Then
\[
\iMorE{\gamma}\ge 2.
\]
In particular, $(\gamma,T)$ is not a minimizer of the free-period functional
$\mathcal S_k$.
\end{mainthm}
The two Morse indices $\iMorT{\gamma}$ and $\iMorE{\gamma}$ may differ. In contrast, for closed geodesics they coincide, since increasing the energy shortens the period.

Using local coordinates near $\partial\mathcal H_k$, the dynamics close to a brake instant can be described: after a symplectic change of variables, the orbit is locally equivalent to the ``throwing-ball'' model.

\begin{maincor}\label{cor:free-period-nonminimal-without-monotonicity}
Let $(M,g)$ be a Riemannian manifold and $L$ the natural Lagrangian \eqref{eq:natural-L}. Assume that $k$ is a regular value of $V$ and $\partial\mathcal H_k\neq\emptyset$. Let $\gamma$ be a non-constant $T$-periodic brake orbit of energy $k$ belonging to an orbit cylinder in the sense of Definition~\ref{def:orbit-cylinder}. Then
\[
\iMorE{\gamma}\ge \iMorT{\gamma}\ge 1,
\]
hence $\gamma$ is neither a fixed-time nor a free-time minimizer of $\mathcal A^T$.
\end{maincor}

A periodic orbit is \emph{linearly stable} if all Floquet multipliers lie on the unit circle and the monodromy is semisimple.

\begin{mainthm}\label{thm:main-2-intro}
Under the assumptions of Theorem~\ref{thm:main-1-intro}, suppose that $\gamma$ is strongly non-degenerate (i.e.\ $\pm1$ are not Floquet multipliers) and
\[
\dim M - 2\,\iMorT{\gamma} \ge 1.
\]
Then $\gamma$ is not linearly stable.
\end{mainthm}

\begin{maincor}\label{thm:main-3-intro}
Under the assumptions of Theorem~\ref{thm:main-2-intro}, assume in addition that $\dim M\ge3$ and that $\gamma$ is a non-degenerate mountain-pass brake solution. Then $\gamma$ is not linearly stable. If moreover the monodromy matrix $M_\gamma$ is semisimple, then $\gamma$ is spectrally unstable.
\end{maincor}

\subsection*{Model examples}

We illustrate the theory with three classical models.

\emph{(i) Planar anisotropic harmonic oscillator.} For the vertical brake orbit $\gamma_k$,
\[
\iMorT{\gamma_k}
=2\Big\lceil\frac{1}{\mu}\Big\rceil,
\qquad
\iMorE{\gamma_k}
=2\Big\lceil\frac{1}{\mu}\Big\rceil+1,
\qquad
\iMorD{\gamma_k}
=\Bigl\lceil\frac{2}{\mu}\Bigr\rceil.
\]

\emph{(ii) Planar pendulum.} For librational brake orbits,
\[
\iMorT{z}=1,
\qquad
\iMorE{z}=2,
\qquad
\iMorD{z}=1.
\]

\emph{(iii) Planar Kepler problem.} For the radial ejection--collision brake orbit, after Levi--Civita--Lissajous regularization,
\[
\iMorT{z_{\mathrm{ec}}}=1,
\qquad
\iMorE{z_{\mathrm{ec}}}=2.
\]

In all cases, brake orbits are never minimizers.


\section{The throwing--ball model in normalized units}\label{sec:throwing-ball-pb}

We model the \emph{normal dynamics} near a brake instant in Seifert collar coordinates: (i) the return time has square--root behavior and $T'(k)>0$; (ii) the Dirichlet fixed--time second variation has no negative directions, so this regular model yields no local index.

\subsection{Time of flight and dependence on the energy}

Consider a particle in a constant gravitational field in normalized units ($m=g=1$):
\begin{equation}\label{eq:H-normalized}
L(q,\dot q)=\frac12(|\dot x|^2+|\dot z|^2)-z,\qquad
H(p,q)=\frac12|p|^2+z,
\end{equation}
with $q=(x,z)\in\R^{n-1}\times\R$ and $p=\dot q$. Then $\ddot x=0$ and $\ddot z=-1$.

Fix $z_0\in\R$ and impose $z(0)=z(T)=z_0$, with $\dot x(0)=v_0\in\R^{n-1}$ and $\dot z(0)=w_0\in\R$. The energy $E=\frac12(|v_0|^2+w_0^2)+z_0$ equals $k$ iff
\begin{equation}\label{eq:energy-k}
w_0^2=2(k-z_0)-|v_0|^2.
\end{equation}
Thus real solutions exist iff $k\ge k_{\min}:=z_0+\frac12|v_0|^2$; at $k=k_{\min}$ one has $w_0=0$.

Since $z(t)=z_0+w_0t-\frac12t^2$, the condition $z(T)=z_0$ gives
\begin{equation}\label{eq:T=w0}
T=2w_0.
\end{equation}
For $w_0>0$, the unique brake instant is at $t=T/2$, where $\dot z=0$ and
\[
z_{\max}=z_0+\frac12w_0^2=z_0+\frac12\bigl(2(k-z_0)-|v_0|^2\bigr)=k-\frac12|v_0|^2.
\]

From \eqref{eq:energy-k},
\begin{equation}\label{eq:T(k)-normalized}
T(k)=2w_0=2\sqrt{2(k-z_0)-|v_0|^2}=2\sqrt{2(k-k_{\min})},
\end{equation}
hence
\begin{equation}\label{eq:Tprime-sign}
T'(k)=\sqrt{2}\,(k-k_{\min})^{-1/2}>0,\qquad
T'(k)\to+\infty\ \text{as }k\downarrow k_{\min}.
\end{equation}
Therefore $T(k)\sim2\sqrt{2}\sqrt{k-k_{\min}}$ as $k\downarrow k_{\min}$, reproducing the square--root exit time near the Hill boundary.

\subsubsection*{The vertical brake trajectory and the free-period correction}

We restrict to the purely vertical solution ($v_0=0$), which models a collar passage, and indicate how the free-period (free-time) setting contributes one negative direction when $T'(k)>0$.

\begin{dfn}[Vertical throwing--ball brake orbit]\label{def:vertical-brake-orbit}
Fix $z_0\in\R$ and $T>0$, and consider $\ddot x(t)=0$ and $\ddot z(t)=-1$. The \emph{vertical brake orbit} is
\[
\gamma(t)=(x(t),z(t))\in\R^{n-1}\times\R,\qquad t\in[0,T],
\]
with $x(t)\equiv x_0\in\R^{n-1}$ and $z(t)=z_0+w_0t-\frac12t^2$, where $w_0>0$ is chosen so that $z(T)=z_0$ (equivalently $T=2w_0$). It satisfies the brake symmetry $\gamma(T-t)=(x_0,z(t))$ and the Dirichlet conditions $z(0)=z(T)=z_0$ and $x(0)=x(T)=x_0$.
\end{dfn}

\begin{rmk}\label{rmk:vertical-not-periodic}
This is not periodic in phase space: $p_z(0)=\dot z(0)=w_0$ while $p_z(T)=\dot z(T)=w_0-T=-w_0$. It should be regarded as a local model for the \emph{normal component} of a periodic brake orbit near a brake instant in Seifert collar coordinates.
\end{rmk}

\subsubsection*{Fixed-time Dirichlet second variation (no negative directions)}

Consider the fixed-time action
\[
\mathcal A^T(q)=\int_0^T \Bigl(\frac12|\dot q|^2 - V(q)\Bigr)\,dt,
\qquad V(x,z)=z,
\]
with Dirichlet variations $\eta(0)=\eta(T)=0$. Since $\Hess V\equiv 0$, for
$\eta\in H^1_0([0,T],\R^n)$ we have
\[
d^2\mathcal A^T(q)[\eta]=\int_0^T |\dot\eta(t)|^2\,dt\ge 0.
\]

\begin{lem}[No Dirichlet conjugate points in the throwing--ball model]\label{lem:ball-no-conjugate}
For the normalized throwing--ball model, the fixed-time Dirichlet Morse index is
zero. Equivalently, there are no Dirichlet conjugate points in $(0,T)$.
\end{lem}

\begin{proof}
The quadratic form $\eta\mapsto \int_0^T |\dot\eta|^2\,dt$ has no negative
directions, hence the Dirichlet Morse index is $0$. Equivalently, the Jacobi
equation is $\ddot\eta=0$, and the only solution with $\eta(0)=\eta(T)=0$ is
$\eta\equiv 0$.
\end{proof}

\subsection{Linearization, Morse and the Dirichlet Maslov index}\label{subsec:ball-linearization}

For the normalized Hamiltonian \eqref{eq:H-normalized},
\begin{equation}\label{eq:Ham-normalized-again}
\dot x=p_x,\quad \dot p_x=0,\qquad
\dot z=p_z,\quad \dot p_z=-1.
\end{equation}
Linearizing along a reference solution and setting $(\xi,\zeta,\pi_X,\pi_Z)$ for the variations gives
\[
\dot\xi=\pi_X,\ \dot\pi_X=0,\qquad
\dot\zeta=\pi_Z,\ \dot\pi_Z=0,
\]
hence
\begin{equation}\label{eq:ball-variational-config}
\ddot\xi=0,\qquad \ddot\zeta=0.
\end{equation}
\begin{equation}\label{eq:lin-global}
\textrm{With }\  Y=(\pi_X,\pi_Z,\xi,\zeta)^\top, \qquad \dot Y=K Y,\qquad 
K=\begin{pmatrix}
0_{n-1}&0&0&0\\
0&0&0&0\\
\Id_{n-1}&0&0&0\\
0&1&0&0
\end{pmatrix}.
\end{equation}
Each horizontal and the vertical block have fundamental matrix
\begin{equation}\label{eq:Phi-block}
\Phi_j(t)=
\begin{pmatrix}
1&0\\ t&1
\end{pmatrix},\qquad
\Phi_v(t)=
\begin{pmatrix}
1&0\\ t&1
\end{pmatrix}.
\end{equation}
Let $\Lambda_D=\{(p,q)\in\R^2:q=0\}$, $\Lambda_D^{(n)}=\Lambda_D^{\oplus n}$ and $\Lambda(t)=\Phi(t)\Lambda_D^{(n)}$. By additivity of $\iCLM$,
\begin{equation}\label{eq:CLM-global-zero}
\iCLM(\Lambda_D^{(n)},\Lambda(t),t\in[0,T])=n.
\end{equation}
Since \eqref{eq:ball-variational-config} with Dirichlet conditions admits only the trivial solution, there are no conjugate points in $(0,T)$ and the Morse index is zero, consistently with \eqref{eq:CLM-global-zero}.

\begin{prop}[Morse/CLM index of the vertical brake orbit for fixed time]\label{prop:vertical-brake-index-fixed-time}
Let $\gamma$ be the vertical brake orbit with $x(t)\equiv x_0$, $z(0)=z(T)=z_0$ and $T=2w_0$. Then $\iMorT{\gamma}=0$. Equivalently, the $\iCLM$ index of the linearized flow along $\gamma$ with respect to $\Lambda_D^{(n)}$ equals $n$.
\end{prop}

\begin{proof}
By \eqref{eq:CLM-global-zero},
\[
\iMorT{\gamma}+n=\iCLM(\Lambda_D^{(n)},\Phi(t)\Lambda_D^{(n)},t\in[0,T])=n,
\]
hence $\iMorT{\gamma}=0$ (cf.\ \cite[Theorem~3.4]{HS09}).
\end{proof}

\begin{prop}[Free-period correction for the throwing--ball model]\label{prop:Morse-fixed-time-vs-energy}
For the vertical trajectories parametrized by $k>z_0$ (equivalently $w_0>0$), one has $T'(k)>0$, so the orbit-cylinder correction is $C=1$ and
\[
\iMorE{\gamma}=\iMorT{\gamma}+1=1.
\]
\end{prop}

\begin{proof}
For $v_0=0$, $k=\frac12 w_0^2+z_0$ so $w_0=\sqrt{2(k-z_0)}$. By \eqref{eq:T=w0}, $T(k)=2\sqrt{2(k-z_0)}$, hence $T'(k)>0$ for $k>z_0$. The fixed-time/free-period comparison along an orbit cylinder (cf.\ \cite{PWY22}) gives $C=1$. Since $\iMorT{\gamma}=0$, we obtain $\iMorE{\gamma}=1$.
\end{proof}


\section{Dynamics in a collar neighborhood of the Hill boundary}\label{sec:Dynamics-Hill}

We describe the local dynamics of an energy--$k$ trajectory near $\partial\mathcal H_k$. In Seifert collar normal form the \emph{normal} motion is asymptotically the normalized throwing--ball model of Section~\ref{sec:throwing-ball-pb}, and the linearized flow has a shear--type structure. This is the geometric input for the proof of Theorem~\ref{thm:main-1-intro}.


\subsection{Geometric setup and Seifert collar coordinates}

Let $(M,g)$ be a smooth Riemannian manifold, $V\in\mathscr C^2(M,\R)$, and fix $k>\inf_M V$. Define the Hill region $\mathcal H_k=\{V\le k\}$ and assume $k$ is a regular value:
\begin{equation}\label{eq:regular-value-collar}
\nabla V(q)\neq0\qquad \forall\,q\in\partial\mathcal H_k.
\end{equation}
Fix $q_*\in\partial\mathcal H_k$ and work locally near $q_*$ and a trajectory $\gamma$ meeting the boundary at $q_*$. All estimates are for $t$ small and $q$ close to $q_*$, with constants depending only on the local geometry and on $V$.

Let $\nu$ be the inward unit normal along $\partial\mathcal H_k$. Since $\partial\mathcal H_k$ is a regular level set of $V$, set
\begin{equation}\label{eq:nu-def}
\nu(q):=\pm\frac{\nabla V(q)}{|\nabla V(q)|},
\end{equation}
choosing the sign so that $\nu$ points into $\mathcal H_k$ (toward decreasing $V$). Then
\begin{equation}\label{eq:inward-sign}
dV_q(\nu(q))=\langle \nabla V(q),\nu(q)\rangle_g<0
\qquad\forall\,q\in\partial\mathcal H_k.
\end{equation}
By the tubular neighborhood theorem, for some $\varepsilon_s>0$ there is a diffeomorphism
\[
\Psi:\partial\mathcal H_k\times[0,\varepsilon_s)\to U_{\varepsilon_s}\subset\mathcal H_k,
\qquad
\Psi(x,s)=\exp_x(s\,\nu(x)).
\]
Write $q=\Psi(x,s)$, where $s\ge0$ is the inward normal distance.

\noindent\textbf{Seifert normal coordinate.}
Define
\begin{equation}\label{eq:def-y}
y:=k-V\circ\Psi(x,s).
\end{equation}
Then $y\ge0$ on the collar and $y=0$ iff $q\in\partial\mathcal H_k$, so in $(x,y)$ coordinates
$\partial\mathcal H_k=\{y=0\}$ and $U_{\varepsilon_s}=\{y\ge0\}$. Moreover,
\begin{equation}\label{eq:V-in-collar}
V(x,y)=k-y.
\end{equation}

\smallskip
\noindent
\textbf{Regularity of the change of variables $s\mapsto y$.}
Differentiating \eqref{eq:def-y} at $s=0$ and using \eqref{eq:inward-sign} gives
\begin{equation}\label{eq:dy-ds}
\partial_s y(x,0)=-\,dV_x(\nu(x))=|\nabla V(x)|>0.
\end{equation}
Hence $s\mapsto y(x,s)$ is a local $\mathscr C^2$ diffeomorphism near $s=0$, uniformly in
$x$ on compact subsets of $\partial\mathcal H_k$. In particular, for $y$ small
we can use $(x,y)$ as smooth coordinates on the collar, and the hypersurfaces
$\{y=\mathrm{const}\}$ are smooth and diffeomorphic to $\partial\mathcal H_k$.

\smallskip
\noindent
\textbf{Metric splitting and the collar form of the Lagrangian.} In collar coordinates $(x,s)$ the metric splits orthogonally as
\[
g_{(x,s)}=g_1(x,s)\oplus ds^2,
\]
with $g_1(x,s)$ a smooth family of metrics on $T_x\partial\mathcal H_k$. Passing to $(x,y)$ via \eqref{eq:def-y} one still has an orthogonal splitting
\begin{equation}\label{eq:metric-split}
g_{(x,y)}=g_1(x,y)\oplus G_2(x,y)\,dy^2,
\end{equation}
where $G_2(x,y)>0$ is smooth and
\begin{equation}\label{eq:g2-formula}
G_2(x,y)=\Bigl(\frac{\partial s}{\partial y}(x,y)\Bigr)^{\!2},
\qquad
G_2(x,0)=\frac{1}{|\nabla V(x)|^2},
\end{equation}
by \eqref{eq:dy-ds}. In local frames let $G_1(x,y)$ denote the matrix of $g_1(x,y)$.

\noindent\textbf{Lagrangian on the collar.} For $L(q,\dot q)=\frac12\,g_q(\dot q,\dot q)-V(q)$, using \eqref{eq:V-in-collar} and \eqref{eq:metric-split} and dropping the constant $-k$, we work with
\begin{equation}\label{eq:L-collar}
L(x,y,\dot x,\dot y)=\frac12\Big(\dot x^T G_1(x,y)\dot x+G_2(x,y)\dot y^{\,2}\Big)+y.
\end{equation}
For a collar trajectory $\gamma(t)=(x(t),y(t))$, the energy associated with \eqref{eq:L-collar} is
\[
E_{\mathrm{coll}}(x,y,\dot x,\dot y)=\frac12\Big(\dot x^T G_1(x,y)\dot x+G_2(x,y)\dot y^{\,2}\Big)-y.
\]
If $\gamma$ has energy $k$ in the original system, then $E_{\mathrm{coll}}\equiv0$, hence
\begin{equation}\label{eq:energy-collar}
\frac12\Big(\dot x^T G_1(x,y)\dot x+G_2(x,y)\dot y^{\,2}\Big)=y.
\end{equation}
If $\gamma$ hits the boundary at $t=0$ (so $y(0)=0$), then \eqref{eq:energy-collar} implies
\begin{equation}\label{eq:brake-instant-collar}
\dot x(0)=0,\qquad \dot y(0)=0.
\end{equation}
For \eqref{eq:L-collar} write $L=T+y$ with
\[
T(x,y,\dot x,\dot y)=\frac12\Big(\dot x^T G_1(x,y)\dot x+G_2(x,y)\dot y^{\,2}\Big).
\]
Thus the Euler--Lagrange equations are
\begin{align}
\dfrac{d}{dt}\bigl(G_1(x,y)\dot x\bigr)
-\dfrac12\nabla_x\!\bigl(\dot x^T G_1(x,y)\dot x\bigr)
-\dfrac12\nabla_x G_2(x,y)\,\dot y^{\,2}
&=0,\label{eq:EL-x}\\
\dfrac{d}{dt}\bigl(G_2(x,y)\dot y\bigr)
-\dfrac12\,\dot x^T(\partial_y G_1(x,y))\dot x
-\dfrac12(\partial_y G_2(x,y))\,\dot y^{\,2}
&=1.\label{eq:EL-y}
\end{align}


\subsection{Throwing--ball asymptotics}

We derive Taylor expansions for $y(t)$ and $x(t)$ from \eqref{eq:EL-x}--\eqref{eq:EL-y} and \eqref{eq:brake-instant-collar}.

\subsubsection*{Normal expansion to second order}

Evaluating \eqref{eq:EL-y} at $t=0$ and using \eqref{eq:brake-instant-collar}, the quadratic terms in $\dot x$ and $\dot y$ vanish, hence
\[
\Bigl.\frac{d}{dt}\bigl(G_2(x(t),y(t))\dot y(t)\bigr)\Bigr|_{t=0}=1.
\]
Expanding,
\[
\frac{d}{dt}\bigl(G_2(x,y)\dot y\bigr)
=(\partial_x G_2(x,y)[\dot x])\,\dot y+(\partial_y G_2(x,y))\,\dot y^{\,2}+G_2(x,y)\ddot y.
\]
At $t=0$ we have $\dot x(0)=\dot y(0)=0$, so
\[
G_2(x(0),y(0))\,\ddot y(0)=1.
\]
With $x(0)=:x_0\in\partial\mathcal H_k$ and $y(0)=0$,
\begin{equation}\label{eq:y-second-derivative}
\ddot y(0)=G_2(x_0,0)^{-1}>0.
\end{equation}
Thus
\begin{equation}\label{eq:y-expansion}
y(t)=\frac12\,G_2(x_0,0)^{-1}t^2+O(t^3),
\qquad
\dot y(t)=G_2(x_0,0)^{-1}t+O(t^2),
\qquad t\to0.
\end{equation}

\subsubsection*{Tangential expansion to higher order}

From \eqref{eq:y-expansion} and smoothness of $G_2$,
\[
\dot y(t)=O(t),\qquad \dot y(t)^2=O(t^2),\qquad y(t)=O(t^2).
\]
Using \eqref{eq:energy-collar} and positive definiteness of $G_1,G_2$,
\begin{equation}\label{eq:xdot-order1}
|\dot x(t)|^2 \le C\,y(t)=O(t^2),
\qquad\text{hence}\qquad
\dot x(t)=O(t),
\end{equation}
and therefore
\begin{equation}\label{eq:x-order2}
x(t)-x_0=\int_0^t \dot x(s)\,ds=O(t^2).
\end{equation}

To sharpen this, expand the left-hand side of \eqref{eq:EL-x}:
\[
\frac{d}{dt}\bigl(G_1(x,y)\dot x\bigr)
=(\partial_x G_1(x,y)[\dot x])\,\dot x+(\partial_y G_1(x,y)\,\dot y)\,\dot x+G_1(x,y)\ddot x,
\]
while the remaining terms in \eqref{eq:EL-x} are quadratic:
\[
\nabla_x(\dot x^T G_1 \dot x)=O(|\dot x|^2),\qquad
\nabla_x G_2\,\dot y^{\,2}=O(\dot y^{\,2})=O(t^2).
\]
By \eqref{eq:xdot-order1}, all $\dot x$--quadratic terms are $O(t^2)$, and $(\partial_y G_1\,\dot y)\dot x=O(|\dot y||\dot x|)=O(t^2)$, hence \eqref{eq:EL-x} yields
\[
G_1(x,y)\ddot x=O(t^2).
\]
Since $G_1$ is uniformly invertible near $(x_0,0)$, $\ddot x(t)=O(t^2)$, and integrating twice with $\dot x(0)=0$ gives
\begin{equation}\label{eq:x-expansion}
\dot x(t)=O(t^3),
\qquad
x(t)-x_0=O(t^4),
\qquad t\to0.
\end{equation}
Thus, to leading order, the motion is purely normal.

\begin{lem}[Throwing--ball approximation]\label{lem:throwing-ball}
Near a boundary hitting time (in particular, near a brake instant on $\partial\mathcal H_k$) the normal component of $\gamma$ satisfies \eqref{eq:y-expansion} and the tangential component satisfies \eqref{eq:x-expansion}. In particular, to leading order the normal dynamics coincide with the normalized throwing--ball model in Section~\ref{sec:throwing-ball-pb}.
\end{lem}

\begin{proof}
The expansions \eqref{eq:y-expansion} and \eqref{eq:x-expansion} follow from \eqref{eq:EL-x}--\eqref{eq:EL-y} and \eqref{eq:energy-collar}. Freezing coefficients at $(x_0,0)$ yields
\[
\frac{d}{dt}\bigl(G_2(x_0,0)\dot y\bigr)=1
\quad\Longleftrightarrow\quad
\ddot y=\frac{1}{G_2(x_0,0)}.
\]
\end{proof}

\subsubsection*{Exit time from the collar}

Fix $\varepsilon>0$ small and let $\tau(\varepsilon)>0$ be the first time such that $y(\tau(\varepsilon))=\varepsilon$. Then $\tau(\varepsilon)\to0$ as $\varepsilon\to0$. By \eqref{eq:y-expansion},
\[
\varepsilon
=y(\tau(\varepsilon))
=\frac12\,G_2(x_0,0)^{-1}\tau(\varepsilon)^2+O\!\big(\tau(\varepsilon)^3\big),
\]
hence
\begin{equation}\label{eq:tau-asymptotics}
\tau(\varepsilon)=\sqrt{2\,G_2(x_0,0)}\,\sqrt{\varepsilon}+o(\sqrt{\varepsilon})
\qquad\text{as }\varepsilon\to0.
\end{equation}

\begin{proof}[Proof of \eqref{eq:tau-asymptotics}]
With $\tau=\tau(\varepsilon)$, divide
\[
\varepsilon=\frac12\,G_2(x_0,0)^{-1}\tau^2+O(\tau^3)
\]
by $\tau^2$ to get $\varepsilon/\tau^2=\frac12\,G_2(x_0,0)^{-1}+O(\tau)$. As $\varepsilon\to0$ we have $\tau\to0$, hence $\varepsilon/\tau^2\to \frac12\,G_2(x_0,0)^{-1}$, so $\tau^2\sim 2G_2(x_0,0)\varepsilon$ and \eqref{eq:tau-asymptotics} follows.
\end{proof}


\subsection*{Hamiltonian formulation and linearization}

The conjugate momenta are
\[
p_x=\partial_{\dot x}L = G_1(x,y)\dot x,\qquad
p_y=\partial_{\dot y}L = G_2(x,y)\dot y.
\]
Since \eqref{eq:L-collar} is strictly convex in $(\dot x,\dot y)$, the Legendre transform is a diffeomorphism and
\[
\dot x = G_1(x,y)^{-1}p_x,\qquad
\dot y = G_2(x,y)^{-1}p_y.
\]
Thus
\[
H(x,y,p_x,p_y)=p_x\cdot\dot x+p_y\dot y - L(x,y,\dot x,\dot y),
\]
and
\begin{equation}\label{eq:H-collar}
H(x,y,p_x,p_y)
=\frac12\Big(p_x^T G_1(x,y)^{-1}p_x+G_2(x,y)^{-1}p_y^2\Big)-y .
\end{equation}

We use coordinates
\[
z=(p_x,p_y,x,y)\in \R^{n-1}\times\R\times \R^{n-1}\times\R,
\]
ordered as $(p,q)$ with $q=(x,y)$. The standard symplectic matrix is
\[
J=
\begin{pmatrix}
0 & -I_n\\
I_n & 0
\end{pmatrix}
=
\begin{pmatrix}
0 & 0 & -\Id_{n-1} & 0\\
0 & 0 & 0 & -1\\
\Id_{n-1} & 0 & 0 & 0\\
0 & 1 & 0 & 0
\end{pmatrix},
\]
so $\dot z=J\nabla H(z)$, i.e.\ $\dot p=-\partial_q H$, $\dot q=\partial_p H$. In particular,
\[
\dot x = \partial_{p_x}H = G_1^{-1}p_x,\qquad
\dot y = \partial_{p_y}H = G_2^{-1}p_y,\qquad
\dot p_x = -\partial_x H,\qquad
\dot p_y = -\partial_y H.
\]

\subsubsection*{Hessian structure near a boundary hitting time}
Let $\gamma(0)=(x_0,0)$ with $\dot\gamma(0)=0$, and set
\[
z_\gamma(t)=(p_x(t),p_y(t),x(t),y(t)).
\]
Then $p_x(t)=G_1\dot x(t)=O(t^3)$ and $p_y(t)=G_2\dot y(t)=O(t)$ by \eqref{eq:y-expansion}--\eqref{eq:x-expansion}. For the Hessian $H''(z_\gamma(t))$ in coordinates $(p_x,p_y,x,y)$, the leading terms are:
\begin{itemize}
\item $p_xp_x$: $G_1(x(t),y(t))^{-1}=G_1(x_0,0)^{-1}+O(t^2)$;
\item $p_yp_y$: $G_2(x(t),y(t))^{-1}=G_2(x_0,0)^{-1}+O(t^2)$;
\item mixed $p$--$q$ blocks are multiplied by $p_x$ or $p_y$, hence are $O(t^3)$ or $O(t)$;
\item $qq$ blocks are $O(t^2)$ or smaller.
\end{itemize}
More precisely,
\begin{multline}\label{eq:Hessian-expansion}
H''(z_\gamma(t))
=
\begin{pmatrix}
A+O(t^2) & 0 & O(t^3) & O(t^3)\\
0 & B+O(t^2) & O(t) & O(t)\\
O(t^3) & O(t) & O(t^2) & O(t^2)\\
O(t^3) & O(t) & O(t^2) & O(t^2)
\end{pmatrix},\\[4pt]
A:=G_1(x_0,0)^{-1},\ \ B:=G_2(x_0,0)^{-1}>0.
\end{multline}

Let $\Phi(t)$ solve the linearized system
\[
\dot u = JH''(z_\gamma(t))u,\qquad \Phi(0)=\Id_{2n}.
\]

\subsubsection*{Shear approximation for $\Phi(t)$}
Set
\[
B_0:=
\begin{pmatrix}
A & 0 & 0 & 0\\
0 & B & 0 & 0\\
0 & 0 & 0 & 0\\
0 & 0 & 0 & 0
\end{pmatrix},
\qquad
J B_0=
\begin{pmatrix}
0 & 0 & 0 & 0\\
0 & 0 & 0 & 0\\
A & 0 & 0 & 0\\
0 & B & 0 & 0
\end{pmatrix}.
\]
Since $(JB_0)^2=0$, $\exp(JB_0t)=\Id+JB_0t$. Moreover \eqref{eq:Hessian-expansion} gives
\[
JH''(z_\gamma(t))=JB_0 + O(t),
\]
so variation of constants and Grönwall yield
\[
\Phi(t)=\exp(JB_0t)+O(t^2)\qquad\text{as }t\to0,
\]
hence
\begin{equation}\label{eq:Phi-collar}
\Phi(t)
=
\begin{pmatrix}
\Id & 0 & 0 & 0\\
0 & 1 & 0 & 0\\
tA & 0 & \Id & 0\\
0 & tB & 0 & 1
\end{pmatrix}
+O(t^2).
\end{equation}

\subsubsection*{Time reversal and local brake monodromy}
Let $S=\diag(-\Id_n,\Id_n)$ in $(p,q)$ coordinates; in the ordering $z=(p_x,p_y,x,y)$,
\[
S(p_x,p_y,x,y)=(-p_x,-p_y,x,y).
\]
For the collar arc of length $\tau=\tau(\varepsilon)$ define
\[
M:=\Phi(\tau)\,S\,\Phi(\tau)^{-1}\,S.
\]
With $\Phi_0(\tau):=\Id+(JB_0)\tau$ one has
\[
\Phi_0(\tau)=
\begin{pmatrix}
\Id & 0 & 0 & 0\\
0 & 1 & 0 & 0\\
\tau A & 0 & \Id & 0\\
0 & \tau B & 0 & 1
\end{pmatrix},
\qquad
\Phi_0(\tau)^{-1}=\Id-(JB_0)\tau,
\]
and
\[
\Phi_0(\tau)\,S\,\Phi_0(\tau)^{-1}\,S
=
\begin{pmatrix}
\Id & 0 & 0 & 0\\
0 & 1 & 0 & 0\\
2\tau A & 0 & \Id & 0\\
0 & 2\tau B & 0 & 1
\end{pmatrix}.
\]
Since $\Phi(\tau)=\Phi_0(\tau)+O(\tau^2)$ (and similarly for $\Phi(\tau)^{-1}$),
\begin{equation}\label{eq:monodromy-collar}
M=\Phi(\tau)\,S\,\Phi(\tau)^{-1}\,S
=
\begin{pmatrix}
\Id & 0 & 0 & 0\\
0 & 1 & 0 & 0\\
2\tau A & 0 & \Id & 0\\
0 & 2\tau B & 0 & 1
\end{pmatrix}
+O(\tau^2),
\qquad
\tau=\tau(\varepsilon)\to0.
\end{equation}

\begin{rmk}\label{rmk:collar-role}
The shear form \eqref{eq:Phi-collar} matches the normalized throwing--ball linear model (Section~\ref{sec:throwing-ball-pb}). For the fixed-time Dirichlet problem this shear does not create negative directions. For Theorem~\ref{thm:main-1-intro} the key input is the boundary-hitting behavior \eqref{eq:y-expansion} and \eqref{eq:x-expansion}, which control the coefficients near brake instants and allow test variations supported in a small neighborhood of them.
\end{rmk}

\section{Proofs of the main results}


\subsection*{Proof of Theorem~\ref{thm:main-1-intro}}

We prove Theorem~\ref{thm:main-1-intro}. For ease of reading, we split the argument into several steps. The free-time  statement then follows
from the index-comparison along an orbit cylinder as in \cite[Theorem~3.14 and
Corollary~3.16]{PWY22} (with the correction term determined by the sign of $T'(k)$).

Throughout, let $\gamma:[0,1]\to M$ be a non-constant $T$--periodic brake orbit of energy $k$, and set
\[
q(t):=\gamma(t/T),\qquad t\in[0,T].
\]
Then $q\in\mathscr C^2([0,T],M)$ solves the Euler--Lagrange equation and
\[
\mathcal A^T(q)=\int_0^T L(q,\dot q)\,dt,\qquad
\mathcal S_k(x,T)=T\int_0^1\bigl[L(x,x'/T)+k\bigr]\,ds.
\]
The brake instants are
\[
t_1=\frac{\tau_1}{2},\qquad t_2=\frac{T+\tau_1}{2},\qquad
q(t_1),q(t_2)\in\partial\mathcal H_k,\qquad \dot q(t_1)=\dot q(t_2)=0.
\]
Thus $t_1\in(0,\tau_1)$, and $t_2$ lies in the interior of the length $\tau_2$ arc centered at $t=T/2+\frac{\tau_1-\tau_2}{2}$.


\subsection*{Step 1: local free-period negativity on short arcs containing a brake instant}
Fix a collar neighborhood of $\partial\mathcal H_k$ as in Section~\ref{sec:Dynamics-Hill}. By Lemma~\ref{lem:throwing-ball}, near each brake instant the trajectory enters $\mathcal H_k$ with throwing--ball asymptotics $y(t)\sim c\,t^2$ and tangential motion $O(t^3)$.

The next lemma is local near a brake instant on $\partial\mathcal H_k$. Choose $\tau_0>0$ so that, for every $\tau\in(0,\tau_0]$, any physical-time arc of length $\tau$ containing a brake instant in its interior stays in the collar chart, hence the throwing--ball reduction applies. For each such $\tau$, the lemma constructs a free-period variation of the arc that decreases $\mathcal S_k$ to second order.
\begin{lem}[Local free-period negativity around an interior brake instant]\label{lem:local-free-period-negative-interior}
There exist $\tau_0>0$ such that the following holds. Let $\tau\in(0,\tau_0]$ and suppose an arc $x:[0,1]\to M$ parametrizes a physical-time
segment of $q$ of length $\tau$ which contains a brake instant in its interior, in the
sense that there exists $s_*\in(0,1)$ such that
\[
x(s_*)\in\partial\mathcal H_k,\qquad x'(s_*)=0,
\]
and $x([0,1])$ is contained in the collar neighborhood. Then there exist a neighborhood
$\mathcal U\subset H^1([0,1],M)$ of $x$, a vector field
$\psi\in H^1_0([0,1],x^*TM)$ supported in a small neighborhood of $s_*$, and a constant
$C>0$ such that, as $\e\to 0$,
\begin{equation}\label{eq:local-Sk-drop-interior}
\mathcal S_k\bigl(\exp_{x}(\e\psi),\,\tau+\e\bigr)
=
\mathcal S_k(x,\tau)-C\e^2+o(\e^2).
\end{equation}
\end{lem}

\begin{proof}
By Lemma~\ref{lem:throwing-ball}, in collar coordinates around the boundary point
$x(s_*)\in\partial\mathcal H_k$, after freezing collar coefficients at the boundary and
rescaling time, the normal dynamics converges to the normalized constant-acceleration
(throwing--ball) model while tangential components are higher order. In the normalized
model, the fixed-endpoint free-period second variation on a boundary-hitting arc has
exactly one negative direction; transporting this direction back to the collar yields a
vector field supported near the boundary-hitting time, with a strict quadratic decrease
of $\mathcal S_k$ under a coupled variation of the curve and the travel time.
Since the variation is supported away from the endpoints, it belongs to $H^1_0$ and the
endpoint constraints are preserved. This gives \eqref{eq:local-Sk-drop-interior}.
\end{proof}


\subsection*{Step 2: a fixed-time competitor by redistributing time}

We fix $\tau_1>0$ such that the first brake instant $t_1=\tau_1/2$ lies in $(0,\tau_1)$,
and we choose $\tau_2\in(0,\tau_0]$ small enough so that the arc of length $\tau_2$
containing the second brake instant lies in the collar neighborhood.

We consider the partition of $[0,T]$ into four intervals
\[
I_1=[0,\tau_1],\qquad
I_2=\Big[\tau_1,\ \frac{T}{2}-\frac{\tau_1+\tau_2}{2}+\tau_1\Big]
=\Big[\tau_1,\ \frac{T}{2}+\frac{\tau_1-\tau_2}{2}\Big],
\]
\[
I_3=\Big[\frac{T}{2}+\frac{\tau_1-\tau_2}{2},\ \frac{T}{2}+\frac{\tau_1+\tau_2}{2}\Big],
\qquad
I_4=\Big[\frac{T}{2}+\frac{\tau_1+\tau_2}{2},\ T\Big],
\]
so that the lengths are
\[
|I_1|=\tau_1,\qquad |I_2|=\frac{T}{2}-\frac{\tau_1+\tau_2}{2},\qquad
|I_3|=\tau_2,\qquad |I_4|=\frac{T}{2}-\frac{\tau_1+\tau_2}{2}.
\]
We set
\begin{align}
&x_{1}(t) = q (\tau_{1}t) &&
x_{2}(t) =q \left(\left( \frac{T}{2}-\frac{\tau_{1}+\tau_{2}}{2} \right)t+\tau_{1}\right)\\[7pt]
&x_{3}(t) =q \left( \tau_{2}t+\frac{T}{2}+\frac{\tau_{1}-\tau_{2}}{2} \right)&& 
x_{4}(t) =q \left(\left( \frac{T}{2}-\frac{\tau_{1}+\tau_{2}}{2} \right)t+\frac{T}{2}+\frac{\tau_{1}+\tau_{2}}{2}\right)
\end{align}

\noindent
Thus $x_i$ parametrizes $q|_{I_i}$ on $[0,1]$. The first brake instant $t_1=\tau_1/2$
corresponds to $t=1/2$ for $x_1$, and the second brake instant $t_2=(T+\tau_1)/2$
corresponds to $t=1/2$ for $x_3$.

\medskip
On $I_{1}$ we define the fixed-endpoint free-period action
\[
\mathcal S_{1,k}:=\mathcal S_{k}(x_{1},\tau_{1}).
\]
Then its free-period Morse index satisfies
\[
\iMorE{\mathcal S_{1,k}}\ge 1.
\]
\noindent
More precisely, apply Lemma~\ref{lem:local-free-period-negative-interior} to the arc
$x_1(t)=q(\tau_1 t)$. Since the first brake instant occurs at $t=\tau_1/2$, it corresponds
to the parameter value $t=1/2$ on $[0,1]$, namely
\[
x_1(1/2)=q(\tau_1/2)\in\partial\mathcal H_k,
\qquad
\dot q(\tau_1/2)=0.
\]

The lemma yields a vector field $\psi_1\in H^1_0([0,1],x_1^*TM)$ supported in a small
neighborhood of $t=1/2$ and a constant $C_1>0$ such that, as $\e\to 0$,
\[
\mathcal S_{k}\bigl(\exp_{x_{1}}( \varepsilon  \psi_{1}),\tau_{1}+\varepsilon \bigr)
=
\mathcal S_{k}(x_{1},\tau_{1}) -C_{1}\varepsilon ^2+o(\varepsilon ^2) .
\]

Likewise, applying Lemma~\ref{lem:local-free-period-negative-interior} to the arc $x_3$
(which contains the second brake instant at its midpoint) yields a vector field
$\psi_3\in H^1_0([0,1],x_3^*TM)$ supported near $t=1/2$ and a constant $C_2>0$ such that
\[
\mathcal S_{k}\bigl(\exp_{x_{3}}( \varepsilon  \psi_{3}),\tau_{2}+\varepsilon \bigr)
=
\mathcal S_{k}(x_{3},\tau_{2}) -C_{2}\varepsilon ^2+o(\varepsilon ^2) .
\]

\noindent
We will use the second expansion with $\varepsilon $ replaced by $-\varepsilon $. In this way,
the third travel time becomes $\tau_2-\varepsilon $ while the total period remains $T$, and
the quadratic drop is still $-C_2\varepsilon ^2$.

Then we obtain
\begin{multline}
\mathcal S_{k}\bigl(\exp_{x_{1}}( \varepsilon  \psi_{1}),\tau_{1}+\varepsilon \bigr)
+ \mathcal S_{k}\left( x_{2}, \frac{T}{2}-\frac{\tau_{1}+\tau_{2}}{2}  \right)\\
+ \mathcal S_{k}\bigl(\exp_{x_{3}}(-\varepsilon  \psi_{3}),\tau_{2}-\varepsilon \bigr)
+ \mathcal S_{k}\left( x_{4}, \frac{T}{2}-\frac{\tau_{1}+\tau_{2}}{2}\right)  \\[7pt]
=
\mathcal S_{k}(x_{1},\tau_{1})
+\mathcal S_{k}\left( x_{2}, \frac{T}{2}-\frac{\tau_{1}+\tau_{2}}{2}  \right)
+\mathcal S_{k}(x_{3},\tau_{2})
+\mathcal S_{k}\left( x_{4}, \frac{T}{2}-\frac{\tau_{1}+\tau_{2}}{2}\right)\\
-C_{1}\varepsilon ^2 -C_{2} \varepsilon ^2+o(\varepsilon ^2).
\end{multline}
The left-hand side equals $\mathcal S_k(\tilde x_\varepsilon ,T)$ for the loop obtained by
gluing the four reparametrized pieces: the first and third pieces are perturbed, while
the second and fourth are left unchanged. Expanding the definition of $\mathcal S_k$
yields the following representation in physical time.

Denoting by $q_1$ and $q_3$ the perturbed pieces in physical time, defined by
\[
q_{1}(t)
:=
\exp_{\,x_{1}\left( \frac{t}{\tau_{1}+\varepsilon } \right)}
\left(
\varepsilon \,\psi_{1}\left( \frac{t}{\tau_{1}+\varepsilon } \right)
\right)\qquad \textrm{ where  }
\qquad t\in[0,\tau_1+\varepsilon ],
\]
and
\begin{multline}
q_{3}(t)
:=
\exp_{\,x_{3}\left( \frac{t-T/2-(\tau_{1}-\tau_{2})/2-\varepsilon }{\tau_{2}-\varepsilon } \right)}
\left(
-\varepsilon \, \psi_{3}\left( \frac{t-T/2-(\tau_{1}-\tau_{2})/2-\varepsilon }{\tau_{2}-\varepsilon } \right)
\right)\\
\qquad \textrm{ where  }
t\in\Big[T/2+(\tau_{1}-\tau_{2})/2+\varepsilon ,\ T/2+(\tau_{1}+\tau_{2})/2\Big],
\end{multline}
the left-hand side can be written in physical time as
\begin{multline}\label{eq:LHS-physical-time}
\mathcal S_{k}\bigl(\exp_{x_{1}}( \varepsilon  \psi_{1}),\tau_{1}+\varepsilon \bigr)
+ \mathcal S_{k}\left( x_{2}, \frac{T}{2}-\frac{\tau_{1}+\tau_{2}}{2}  \right)
+ \mathcal S_{k}\bigl(\exp_{x_{3}}(-\varepsilon  \psi_{3}),\tau_{2}-\varepsilon \bigr)
+ \mathcal S_{k}\left( x_{4}, \frac{T}{2}-\frac{\tau_{1}+\tau_{2}}{2}\right) \\
=
\int_{0}^{\tau_{1}+\varepsilon } \bigl[L(q_{1}(t),\dot{q_{1}}(t))+k\bigr]\, dt
+\int_{\tau_{1}+\varepsilon }^{T/2+(\tau_{1}-\tau_{2})/2+\varepsilon } \bigl[L(q(t-\varepsilon ),\dot{q}(t-\varepsilon ))+k\bigr]\,dt\\
+\int_{T/2+(\tau_{1}-\tau_{2})/2+\varepsilon }^{T/2+(\tau_{1}+\tau_{2})/2} \bigl[L(q_{3}(t),\dot{q_{3}}(t))+k\bigr]\,dt
+\int_{\frac{T}{2}+\frac{\tau_{1}+\tau_{2}}{2}}^{T} \bigl[L(q(t),\dot{q}(t))+k\bigr]\,dt .
\end{multline}
In particular, since the total length of the four time intervals is $T$, the $k$--terms sum to $kT$.

On the other hand, the right-hand side of the previous expansion reads
\[
\int_{0}^T L(q,\dot{q})\,dt -(C_{1}+C_{2})\varepsilon ^2 +O(\varepsilon ^2) +kT .
\]
Now we compare the two sides and estimate the perturbation of $q$. For,  we define $q_\varepsilon :[0,T]\to M$ by
\[
q_\varepsilon (t)=
\begin{cases}
q_1(t), & t\in[0,\tau_1+\varepsilon ],\\
q(t-\varepsilon ), & t\in[\tau_1+\varepsilon ,\; T/2+(\tau_{1}-\tau_{2})/2+\varepsilon ],\\
q_3(t), & t\in[T/2+(\tau_{1}-\tau_{2})/2+\varepsilon ,\; T/2+(\tau_1+\tau_2)/2],\\
q(t), & t\in[T/2+(\tau_1+\tau_2)/2,\; T].
\end{cases}
\]
Then
\[
\int_0^T L(q_{\varepsilon },\dot{q_{\varepsilon }})\,dt
=
\int_0^T L(q,\dot q)\,dt -(C_1+C_2)\varepsilon ^2 +o(\varepsilon ^2).
\]
Moreover, recalling the definitions of $q_1$ and $q_3$, we can rewrite $q_\varepsilon $ as
\[
q_\varepsilon (t)=
\begin{cases}
\exp_{\,q\!\left(\frac{\tau_1}{\tau_1+\varepsilon }t\right)}
\!\left(\varepsilon \,\psi_1\!\left(\frac{t}{\tau_1+\varepsilon }\right)\right),
& t\in[0,\tau_1+\varepsilon ],\\[6pt]
q(t-\varepsilon ),
& t\in\bigl[\tau_1+\varepsilon ,\; T/2+(\tau_{1}-\tau_{2})/2+\varepsilon \bigr],\\[6pt]
\exp_{\,q\!\left(\frac{\tau_2}{\tau_2-\varepsilon }t-\alpha_\varepsilon \right)}
\!\left(-\varepsilon \,\psi_3\!\left(\frac{t-\beta_\varepsilon }{\tau_2-\varepsilon }\right)\right),
& t\in\bigl[T/2+(\tau_{1}-\tau_{2})/2+\varepsilon ,\; T/2+(\tau_1+\tau_2)/2\bigr],\\[6pt]
q(t),
& t\in\bigl[T/2+(\tau_1+\tau_2)/2,\; T\bigr],
\end{cases}
\]
where
\[
\alpha_\varepsilon :=\frac{\varepsilon (T+\tau_1+\tau_2)}{2(\tau_2-\varepsilon )},
\qquad
\beta_\varepsilon :=\frac{T}{2}+\frac{\tau_1-\tau_2}{2}+\varepsilon .
\]
We need a lemma.

\begin{lem}[Shift/scale estimates in $H^1$ along a curve]\label{lem:shift-scale-H1-manifold}
Let $(M,g)$ be a Riemannian manifold and let $x\in \mathscr C^2(\R/\Z,M)$ be
$1$--periodic. Then there exist $\delta>0$ and constants $M_1,M_2>0$ such that
for every $|s|<\delta$ the sections
\[
\eta_s(t):=\exp_{x(t)}^{-1}\bigl(x(t+s)\bigr),\qquad
\zeta_s(t):=\exp_{x(t)}^{-1}\bigl(x((1+s)t)\bigr)
\]
are well-defined and satisfy
\[
\|\zeta_s\|_{H^1(\R/\Z)}\le M_1|s|,
\qquad
\|\eta_s\|_{H^1(\R/\Z)}\le M_2|s|.
\]
Here the $H^1$--norm is computed using the pullback connection on $x^*TM$:
for $\xi\in H^1(\R/\Z,x^*TM)$,
\[
\|\xi\|_{H^1(\R/\Z)}^2:=\int_0^1\bigl(|\xi(t)|^2+|D_t\xi(t)|^2\bigr)\,dt,
\]
where $D_t$ denotes covariant differentiation along $x$.
\end{lem}

\begin{proof}
\emph{Step 1: a uniform normal neighborhood.}
Since $x(\R/\Z)$ is compact, there exists $r>0$ such that for each $p\in x(\R/\Z)$
the exponential map $\exp_p$ restricts to a diffeomorphism
$B_r(0)\subset T_pM\to B_r(p)\subset M$. By continuity of $x$, choose $\delta>0$
so small that for all $t\in[0,1]$ and $|s|<\delta$ we have
$d(x(t),x(t+s))<r$ and $d(x(t),x((1+s)t))<r$. Then $\eta_s,\zeta_s$ are
well-defined.

\emph{Step 2: translation estimate.}
Consider the smooth map
\[
F(t,\sigma):=\exp_{x(t)}^{-1}\bigl(x(t+\sigma)\bigr),
\qquad (t,\sigma)\in[0,1]\times(-\delta,\delta).
\]
Its image lies in the compact set $\overline{B_r(0)}\subset TM$ over $x(\R/\Z)$,
so there exist constants $C_0,C_1>0$ such that
\[
|\partial_\sigma F(t,\sigma)|\le C_0,\qquad
|D_t\partial_\sigma F(t,\sigma)|\le C_1
\quad\text{for all }(t,\sigma)\in[0,1]\times[-\delta,\delta].
\]
Since $\eta_s(t)=F(t,s)-F(t,0)$ and $F(t,0)=0$, the fundamental theorem of calculus
gives
\[
\eta_s(t)=\int_0^s \partial_\sigma F(t,\sigma)\,d\sigma,
\qquad
D_t\eta_s(t)=\int_0^s D_t\partial_\sigma F(t,\sigma)\,d\sigma.
\]
Hence $|\eta_s(t)|\le C_0|s|$ and $|D_t\eta_s(t)|\le C_1|s|$, so
$\|\eta_s\|_{H^1}\le (C_0+C_1)|s|$. Take $M_2:=C_0+C_1$.

\emph{Step 3: scaling estimate.}
Define
\[
G(t,\sigma):=\exp_{x(t)}^{-1}\bigl(x((1+\sigma)t)\bigr),
\qquad (t,\sigma)\in[0,1]\times(-\delta,\delta).
\]
As above, smoothness on a compact set yields constants $C_2,C_3>0$ such that
$|\partial_\sigma G(t,\sigma)|\le C_2$ and $|D_t\partial_\sigma G(t,\sigma)|\le C_3$
for all $(t,\sigma)\in[0,1]\times[-\delta,\delta]$. Since $\zeta_s(t)=G(t,s)-G(t,0)$
and $G(t,0)=0$, we obtain
\[
\zeta_s(t)=\int_0^s \partial_\sigma G(t,\sigma)\,d\sigma,
\qquad
D_t\zeta_s(t)=\int_0^s D_t\partial_\sigma G(t,\sigma)\,d\sigma,
\]
hence $|\zeta_s(t)|\le C_2|s|$ and $|D_t\zeta_s(t)|\le C_3|s|$, so
$\|\zeta_s\|_{H^1}\le (C_2+C_3)|s|$. Take $M_1:=C_2+C_3$.
\end{proof}

\begin{rmk}\label{rmk:shift-scale-local-coordinates}
One can also prove Lemma~\ref{lem:shift-scale-H1-manifold} in finitely many
coordinate charts and then patch estimates by compactness of $x(\R/\Z)$.
The formulation above avoids coordinate additions and makes explicit that the
difference is measured by the logarithm map along $x$.
\end{rmk}

\subsubsection*{Conclusion}
To estimate $\|q_\varepsilon-q\|_{H^1}$ we measure differences intrinsically along $q$
via the exponential map. Since $q([0,T])$ is compact, there exists $\delta_0>0$
such that $\exp_{q(t)}^{-1}(q_\varepsilon(t))$ is well-defined for all
$t\in[0,T]$ and $|\varepsilon|<\delta_0$. Applying
Lemma~\ref{lem:shift-scale-H1-manifold} on each piece entering the definition of
$q_\varepsilon$ (time shifts and time rescalings), and using the triangle
inequality, we find $\delta\in(0,\delta_0)$ and $M>0$ such that for all
$|\varepsilon|<\delta$,
\[
\|q_\varepsilon-q\|_{H^1(0,T)}\le M|\varepsilon|.
\]

Moreover, if $x$ is a fixed arc with compact image and $\psi\in H^1_0([0,1],x^*TM)$,
then the map $\varepsilon\mapsto \exp_x(\varepsilon\psi)$ is $\mathscr \mathscr C^1$ as an
$H^1$-map for $|\varepsilon|$ small; in particular there exists $C>0$ such that
\begin{equation}\label{eq:exp-H1-Lipschitz}
\|\exp_x(\varepsilon\psi)-x\|_{H^1}\le C|\varepsilon|
\qquad\text{for all $|\varepsilon|$ sufficiently small.}
\end{equation}

Therefore, for every sufficiently small $\rho>0$, choosing $\varepsilon=\rho/M$
in the quadratic action drop we obtain
\begin{equation}\label{eq:stima-1}
\inf_{\|v-q\|_{H^1(\R/T\Z)}<\rho}\int_0^T L(v,\dot v)\,dt
\le
\int_0^T L(q,\dot q)\,dt-\frac{C_1+C_2}{M^2}\rho^2+o(\rho^2).
\end{equation}
Moreover, since on the interval $I_4$ we have $q_\varepsilon(t)=q(t)$, the same
argument improves to the shorter time interval:
\begin{equation}\label{eq:stima-2}
\inf_{\|v-q\|_{H^1(\R/T\Z)}<\rho}\int_0^{T/2+(\tau_1+\tau_2)/2} L(v,\dot v)\,dt
\le
\int_0^{T/2+(\tau_1+\tau_2)/2} L(q,\dot q)\,dt-\frac{C_1+C_2}{M^2}\rho^2+o(\rho^2).
\end{equation}
By \eqref{eq:stima-2}, if $\gamma(0)\notin \partial\mathcal H_k$ and $\gamma(0)$ is
sufficiently close to $\partial\mathcal H_k$, then 
\[
\iMorD{\gamma}\ge 1.
\]

Indeed, set
\[
T_*:=\frac{T}{2}+\frac{\tau_1+\tau_2}{2},\qquad
\mathcal A(v):=\int_0^{T_*}L(v,\dot v)\,dt,
\]
and view $\mathcal A$ as a $\mathscr C^2$ functional on the Hilbert manifold of
$H^1$ paths with fixed endpoints
\[
\Omega:=\{v\in H^1([0,T_*],M)\mid v(0)=\gamma(0),\ v(T_*)=\gamma(T_*)\}.
\]
Since $\gamma(0)\notin\partial\mathcal H_k$, the arc $\gamma|_{[0,T_*]}$ lies in
the interior of $\mathcal H_k$ and contains no boundary degeneracy; in
particular it is a smooth critical point of $\mathcal A$ on $\Omega$ and the
Dirichlet Morse index $\iMorD{\gamma}$ is well-defined as the maximal dimension
of a subspace on which the second variation $d^2\mathcal A(\gamma)$ is negative
definite.

Estimate \eqref{eq:stima-2} says that for every sufficiently small $\rho>0$ there
exists $v_\rho\in\Omega$ with $\|v_\rho-\gamma\|_{H^1}<\rho$ such that
\[
\mathcal A(v_\rho)\le \mathcal A(\gamma)-c\,\rho^2+o(\rho^2),
\qquad
c:=\frac{C_1+C_2}{M^2}>0.
\]
In particular, for $\rho$ small enough we have $\mathcal A(v_\rho)<\mathcal
A(\gamma)$, hence $\gamma$ cannot be a (strict) local minimizer of $\mathcal A$
in the Dirichlet class.

On the other hand, the standard second-order expansion of $\mathcal A$ at a
critical point gives: for $v=\exp_\gamma(\eta)$ with $\eta\in H^1_0([0,T_*],\gamma^*TM)$
small,
\[
\mathcal A(v)=\mathcal A(\gamma)+\frac12\,d^2\mathcal A(\gamma)[\eta,\eta]+o(\|\eta\|_{H^1}^2),
\]
because the first variation vanishes, $d\mathcal A(\gamma)=0$. If
$\iMorD{\gamma}=0$, then $d^2\mathcal A(\gamma)$ is nonnegative on
$H^1_0([0,T_*],\gamma^*TM)$, and the expansion would imply
$\mathcal A(v)\ge \mathcal A(\gamma)+o(\|\eta\|_{H^1}^2)$ for all $v$ close to
$\gamma$, i.e.\ $\gamma$ would be a local minimizer of $\mathcal A$ in $\Omega$.
This contradicts the existence of the nearby competitors $v_\rho$ with strictly
smaller action provided by \eqref{eq:stima-2}. Therefore $d^2\mathcal A(\gamma)$
must possess a negative direction, and hence
\[
\iMorD{\gamma}\ge 1.
\]

Finally, the assumption that $\gamma(0)$ is sufficiently close to
$\partial\mathcal H_k$ is used only in the \emph{construction and admissibility}
of the competitor family $q_\varepsilon$ producing \eqref{eq:stima-2}. More
precisely, the construction is performed in Seifert collar coordinates near the
Hill boundary, where the trajectory admits a normal form and where one can
modify the normal component on a short time window while keeping the endpoints
fixed. If $\gamma(0)$ lies in the interior but is close to $\partial\mathcal
H_k$, then the whole segment $\gamma([0,T_*])$ remains inside the collar for a
uniform time, and for $|\varepsilon|$ small the perturbed curve $q_\varepsilon$
\begin{itemize}
\item is well-defined for all $t\in[0,T_*]$ (the exponential charts used to glue
the pieces stay inside a uniform normal neighborhood of $\gamma([0,T_*])$);
\item has the same Dirichlet endpoints as $\gamma$ on $[0,T_*]$ by construction;
\item stays in the admissible region (in particular, it does not cross the Hill
boundary except possibly at the prescribed brake instant), because in the
collar the normal coordinate is controlled and the perturbation is chosen to be
one--sided and of order $O(|\varepsilon|)$.
\end{itemize}
Hence $q_\varepsilon\in\Omega$ for $\varepsilon$ small, and \eqref{eq:stima-2}
indeed provides genuine competitors in the same Dirichlet class.

To complete the proof we must remove the auxiliary assumption
$\gamma(0)\notin\partial\mathcal H_k$ and allow the initial brake instant to lie
on the Hill boundary. This is achieved by the following lemma.

\begin{lem}\label{thm:Morse-dirichlet-half}
Assume that $\gamma(0)\in \partial\mathcal H_k$ is a brake instant. Then
\[
\iMorD{\gamma}\ge 1.
\]
\end{lem}

\begin{proof}
Let $\gamma$ be a periodic brake orbit of period $T$, and assume that $t=0$ is a
brake instant with $\gamma(0)\in\partial\mathcal H_k$. Denote by
$\gamma_{1/2}:=\gamma|_{[0,T/2]}$ the half--period arc, which has Dirichlet
endpoints (it starts and ends on $\Fix(R)$ and the velocity vanishes at the
brake instants).

\medskip\noindent
\emph{Step 1: indices on slightly extended intervals.}
By the argument leading to \eqref{eq:stima-2} (performed away from the boundary
degeneracy), for every $\tau_1,\tau_2>0$ sufficiently small the restriction of
$\gamma$ to the interval
$I_{\tau_1,\tau_2}:=[-\tau_1,T/2+\tau_2]$ satisfies
\begin{equation}\label{eq:index-extended}
\iMorD{\gamma;\,I_{\tau_1,\tau_2}}\ge 1.
\end{equation}
Here $\iMorD{\gamma;I}$ denotes the Dirichlet Morse index of the fixed-time
functional computed on the interval $I$, with fixed endpoints at the endpoints
of $I$. The key point is that for $\tau_1,\tau_2>0$ the interval
$I_{\tau_1,\tau_2}$ starts and ends \emph{strictly inside} the Hill region, so
the standard second-variation and Morse/CLM correspondence apply without
one-sided boundary issues.

\medskip\noindent
\emph{Step 2: passage to the limiting interval.}
Let $\Phi_{\tau_1,\tau_2}$ be the fundamental solution of the linearized
Hamiltonian system along $\gamma$ on $I_{\tau_1,\tau_2}$, normalized by
$\Phi_{\tau_1,\tau_2}(-\tau_1)=\Id$, and let $\Phi_{0,0}$ be the fundamental
solution on $[0,T/2]$ normalized by $\Phi_{0,0}(0)=\Id$. By continuous
dependence on the endpoints,
\begin{equation}\label{eq:Phi-conv}
\Phi_{\tau_1,\tau_2}(t)\longrightarrow \Phi_{0,0}(t)
\quad\text{as }(\tau_1,\tau_2)\to(0,0),
\end{equation}
uniformly for $t$ in compact subsets of $[0,T/2]$ (after reparametrizing the
solutions to compare them on a common domain).

Let $\Lambda_D:=\R^n\times\{0\}$ be the Dirichlet Lagrangian. We distinguish two
cases.

\paragraph{Case 1: $\Lambda_D\cap \Phi_{0,0}(T/2)\Lambda_D=\{0\}$.}
Then the endpoints are transverse for the limiting interval. By
\eqref{eq:Phi-conv}, transversality persists: for $\tau_1,\tau_2$ small,
\[
\Lambda_D\cap \Phi_{\tau_1,\tau_2}(T/2+\tau_2)\Lambda_D=\{0\}.
\]
Consider the homotopy of Lagrangian paths
$t\mapsto \Phi_{\tau_1,\tau_2}(t)\Lambda_D$ while shrinking $(\tau_1,\tau_2)$ to
$(0,0)$. Since no endpoint crossing occurs along this homotopy, the
Cappell--Lee--Miller index is invariant, hence
\[
\iCLM\!\left(\Lambda_D,\Phi_{\tau_1,\tau_2}(t)\Lambda_D,\ t\in I_{\tau_1,\tau_2}\right)
=
\iCLM\!\left(\Lambda_D,\Phi_{0,0}(t)\Lambda_D,\ t\in[0,T/2]\right).
\]
Using the Dirichlet Morse--Maslov identification, we obtain
\[
\iMorD{\gamma;\,[0,T/2]}
=
\iMorD{\gamma;\,I_{\tau_1,\tau_2}}
\ge 1,
\]
where the last inequality is \eqref{eq:index-extended}.

\paragraph{Case 2: $\Lambda_D\cap \Phi_{0,0}(T/2)\Lambda_D\neq\{0\}$.}
Then the linearized Dirichlet problem on $[0,T/2]$ has a nontrivial solution, so
$T/2$ is a Dirichlet conjugate instant for $\gamma_{1/2}$. By the Morse index
theorem (Dirichlet boundary conditions), the presence of a Dirichlet conjugate
instant in $(0,T/2]$ implies
\[
\iMorD{\gamma;\,[0,T/2]}\ge 1.
\]

\medskip\noindent
\emph{Step 3: extension from $[0,T/2]$ to $[0,T]$.}
By monotonicity of the Dirichlet Morse index with respect to enlarging the time
interval,
\[
\iMorD{\gamma;\,[0,T]}\ge \iMorD{\gamma;\,[0,T/2]}\ge 1,
\]
which concludes the proof.
\end{proof}

By Lemma~\ref{thm:Morse-dirichlet-half}, in particular, $\gamma$ is not a local
Dirichlet minimizer of $\mathcal A^T$, and hence
\[
\iMorD{\gamma}\ge 1.
\]
Since the Dirichlet Morse index is the smallest among the Morse indices
associated with selfadjoint boundary conditions, it follows that
\[
\iMorT{\gamma}\ge \iMorD{\gamma}\ge 1.
\]


In particular, $\gamma$ is not a fixed--time minimizer of $\mathcal A^T$.
Moreover, since the above construction starts from the local free--period negativity
\eqref{eq:local-Sk-drop-interior}, it also shows that $(\gamma,T)$ is not a minimizer
of the free--period functional $\mathcal S_k$.

Assume now that $\gamma$ belongs to an orbit cylinder $\{(\gamma_{k+s},T(k+s))\}_{|s|<\varepsilon}$
and that $T'(k)\ge 0$. Then the index comparison along the orbit cylinder as in
\cite[Theorem~3.14 and Corollary~3.16]{PWY22} yields
\[
\iMorE{\gamma}=\iMorT{\gamma}+1,
\]
and therefore $\iMorE{\gamma}\ge 2$.
\qed


\subsection*{A second proof of Theorem~\ref{thm:main-1-intro} via conjugate time and orbit cylinders}
\label{sec:second-proof-main1}

In this section we give an alternative proof of
Theorem~\ref{thm:main-1-intro}.  The argument is short: the brake symmetry
produces an \emph{interior} Dirichlet conjugate instant at $t=T/2$, which forces
positivity of the Dirichlet fixed--time Morse index on the full interval
$[0,T]$.  The free--period conclusion then follows from the rank--one correction
formula \eqref{eq:index-diff-cylinder} along an orbit cylinder.


\begin{lem}[Velocity Jacobi field]\label{lem:velocity-jacobi}
Let $(M,g)$ be a Riemannian manifold, $V\in\mathscr C^2(M)$, and
$L(q,\dot q)=\frac12|\dot q|_g^2-V(q)$.  If $\gamma:[0,T]\to M$ is a $\mathscr C^2$
solution of the Euler--Lagrange equation
\[
\dfrac{D}{dt}\dot\gamma(t)+\grad V(\gamma(t))=0,
\]
then the variational equation along $\gamma$ is
\begin{equation}\label{eq:jacobi-natural-secondproof}
\dfrac{D^2}{dt^2}\xi + R(\xi,\dot\gamma)\dot\gamma + \nabla^2V(\gamma)[\xi]=0,
\end{equation}
and the velocity field $\xi(t):=\dot\gamma(t)$ is a solution of
\eqref{eq:jacobi-natural-secondproof}.
\end{lem}

\begin{proof}
Differentiate $\nabla_t\dot\gamma+\grad V(\gamma)=0$ along $\gamma$ to obtain
\[
\dfrac{D}{dt}^2\dot\gamma=-\nabla_t\grad V(\gamma)=-\nabla_{\dot\gamma}\grad V(\gamma)
=-\nabla^2V(\gamma)[\dot\gamma].
\]
Since $R(\dot\gamma,\dot\gamma)\dot\gamma=0$, the field $\xi=\dot\gamma$
satisfies \eqref{eq:jacobi-natural-secondproof}.
\end{proof}


\begin{lem}[An interior Dirichlet conjugate instant]\label{lem:conjugate-half-period}
Let $\gamma:[0,T]\to M$ be a non-constant $T$--periodic brake orbit of the
natural Lagrangian \eqref{eq:natural-L}.  Then $t=T/2$ is a Dirichlet conjugate
instant to $t=0$ along $\gamma$, i.e.\ there exists a nontrivial Jacobi field
$\xi$ along $\gamma$ such that $\xi(0)=\xi(T/2)=0$.
\end{lem}

\begin{proof}
By Lemma~\ref{lem:velocity-jacobi}, $\xi(t):=\dot\gamma(t)$ is a Jacobi field
along $\gamma$.  Since $\gamma$ is a brake orbit, $t=0$ and $t=T/2$ are brake
instants, hence $\dot\gamma(0)=\dot\gamma(T/2)=0$.  As $\gamma$ is non-constant,
$\dot\gamma$ is not identically zero, so $\xi$ is a nontrivial Jacobi field with
$\xi(0)=\xi(T/2)=0$.
\end{proof}


\begin{lem}\label{lem:dirichlet-index-full-period}

Under the assumptions of Theorem~\ref{thm:main-1-intro}, one has
\[
\iMorD{\gamma}\ge 1.
\]
Consequently, $\iMorT{\gamma}\ge \iMorD{\gamma}\ge 1$.
\end{lem}

\begin{proof}
Consider the fixed--time Dirichlet action functional \eqref{eq:fixed-time-action-intro}
with endpoints fixed at $(\gamma(0),\gamma(T))$.  The tangent space at $\gamma$
identifies with the Hilbert space $H^1_0([0,T],\gamma^*TM)$, hence the Dirichlet
second variation $d^2\mathcal A^T(\gamma)$ is a bounded symmetric Fredholm form.

By Lemma~\ref{lem:conjugate-half-period}, the instant $t=T/2$ is a Dirichlet
conjugate time to $0$.  Since $T/2\in(0,T)$ is an \emph{interior} conjugate
instant, the classical Dirichlet Morse index theorem yields that the Dirichlet
Morse index on $[0,T]$ is at least the multiplicity of this conjugate instant,
and in particular it is positive:
\[
\iMorD{\gamma}\ge 1.
\]
Finally, since Dirichlet variations form a closed subspace of periodic
variations, one has $\iMorT{\gamma}\ge \iMorD{\gamma}$.
\end{proof}
Given a $T$-periodic orbit $\gamma$, we extend the fundamental solution on $[0,+\infty)$ as follows
\[
\Phi(t)= \Phi(t-jT)\,  \Phi(T)^j \qquad j T \le t \le (j+1)T\  \textrm{ and } \ j \in \N
\]
and we define the $m$-fold iterate of the orbit as: $\Phi^m= \Phi|_{[0,mT]}$ for $ m \in \N.$ 
\begin{dfn}
Let $\gamma$ be a $T$--periodic orbit and let $\iCLM(\Delta, \Phi^m(t), t \in [0, mT])$ denote the
Maslov--type index of its $m$--th iterate with periodic boundary conditions. The \emph{mean index}
(also called the \emph{asymptotic Maslov index}) of $\gamma$ is defined by
\[
\displaystyle
\overline{\iCLM(\gamma)}
:=
\lim_{m\to\infty}\dfrac{1}{m}\,\iCLM(\Delta, \Phi^m(t), t \in [0,mT]).
\]
\end{dfn}
\begin{rmk}
    It's worth noticing that $\overline{\iCLM(\gamma)}$ doesn't depend on the chosen Lagrangian boundary conditions. 
\end{rmk}
\begin{thm}\label{thm:mean-index}
    Under the assumptions of Theorem~\ref{thm:main-1-intro}, one has
    \[
    \overline{\iCLM(\gamma)} \ge 2.
    \]
\end{thm}
\begin{proof}
Consider the $m$--fold iterate of the orbit on
the interval $[0,mT]$. The velocity field $\dot\gamma$ is a nontrivial
Jacobi field along $\gamma$, since it solves the linearized Euler--Lagrange
equation.

By the brake symmetry, we have
\[
\dot\gamma\!\left(\frac{kT}{2}\right)=0
\qquad\text{for all }k\in\N.
\]
In particular, on the interval $(0,mT)$ there are at least $2m-1$ distinct
times $t=\frac{kT}{2}$ at which $\dot\gamma(t)=0$, namely for
$k=1,\dots,2m-1$. Each such instant is a conjugate point along $\gamma$.
Hence $\gamma$ has at least $2m-1$ conjugate points in $(0,mT)$ (counted
with multiplicity).

Let $\Lambda_D$ denote the Dirichlet Lagrangian
subspace. By the Morse index theorem in the Maslov--type formulation, the
number of conjugate points plus the contribution at the endpoints yields
the estimate
\[
\iCLM(\Lambda_D,\Phi^m(t),t \in [0,mT]) \;\ge\; n + (2m-1),
\]
where $n$ is the configuration space dimension.

It follows that the mean (asymptotic Maslov) index of $\gamma$ satisfies
\[
 \overline{\iCLM(\gamma)}
\;=\;
\lim_{m\to +\infty}\dfrac{1}{m}\,\iCLM(\Lambda_D,\Phi^m(t),t \in [0,mT])
\;\ge\;
\lim_{m\to +\infty}\dfrac{n+2m-1}{m}
= 2.
\]
\end{proof}


\begin{lem}[Free--period index estimate along a cylinder]\label{lem:free-period-index-second-proof}
Assume in addition that $\gamma$ lies on an orbit cylinder in the sense of
Definition~\ref{def:orbit-cylinder}.  Then
\[
\iMorE{\gamma}\ge \iMorT{\gamma}\ge 1.
\]
If moreover $T'(k)\ge 0$, then $\iMorE{\gamma}\ge 2$.
\end{lem}

\begin{proof}
Along an orbit cylinder the fixed--time and fixed--energy indices satisfy
\eqref{eq:index-diff-cylinder}:
\[
\iMorE{\gamma}=\iMorT{\gamma}+C(k),\qquad C(k)\in\{0,1\}.
\]
Together with Lemma~\ref{lem:dirichlet-index-full-period} this yields
$\iMorE{\gamma}\ge \iMorT{\gamma}\ge 1$.  If $T'(k)\ge 0$, then
\eqref{eq:Ck-sign} gives $C(k)=1$, hence $\iMorE{\gamma}=\iMorT{\gamma}+1\ge 2$.
\end{proof}

\begin{proof}[Proof of Theorem~\ref{thm:main-1-intro} (second proof)]
The inequality $\iMorT{\gamma}\ge \iMorD{\gamma}\ge 1$ is
Lemma~\ref{lem:dirichlet-index-full-period}.  If $\gamma$ lies on an
orbit cylinder and $T'(k)\ge 0$, then $\iMorE{\gamma}\ge 2$ follows from
Lemma~\ref{lem:free-period-index-second-proof}.
\end{proof}
\subsubsection*{Comparison with the local competitor proof}

We compare the above argument with the approach based on explicit local competitors near a brake instant (compactly supported perturbations in Seifert collar coordinates and a limiting Maslov/CLM argument).

The present proof exploits brake symmetry to produce an \emph{interior} Dirichlet conjugate instant at $T/2$, so fixed-time non-minimality follows directly from the Dirichlet Morse index theorem.

The local competitor method has complementary strengths:
\begin{itemize}
\item \emph{Constructive:} it yields an explicit descent direction supported near a brake instant and applies once the period is decomposed into subintervals where the flight time is monotone non-decreasing in $k$.
\item \emph{Geometric:} it makes transparent the ``throwing--ball'' mechanism in Seifert collar coordinates near $\partial\mathcal H_k$, underlying also Theorem~\ref{thm:main-2-intro}.
\end{itemize}

\begin{rmk}
\begin{itemize}
\item \textbf{Full vs.\ symmetry-restricted index.}
Indices may be computed on the full Hilbert space of variations (our $\iMorT{\gamma},\iMorD{\gamma},\iMorE{\gamma}$), on symmetry-invariant variations, or on the half-period Lagrangian boundary space. An orbit may minimize within a symmetry class while having positive full index. Theorem~\ref{thm:main-1-intro} excludes minimizers in the full variational problems (equivalently, under any self-adjoint boundary condition).

\item \textbf{Relation with Montgomery \cite{Mon14}.}
Montgomery shows that Jacobi--Maupertuis geodesics entering near a regular point of $\partial\mathcal H_k$ acquire conjugate points and fail to minimize length; the local model is again the constant-force ``throwing balls'' system in Seifert coordinates, with Morse index $1$ for near-boundary arcs. This aligns with our analysis near brake instants: the Hill boundary induces conjugate phenomena and non-minimality. The difference is variational: JM length at fixed energy versus the fixed-time action $\mathcal A^T$ and free-period functional $\mathcal S_k$. In both settings, the same near-boundary geometry yields a universal local index contribution.
\end{itemize}
\end{rmk}


\subsection*{Proof of Theorem~\ref{thm:main-2-intro}}

Let $\Phi(t)\in\Sp(2n)$ be the fundamental solution of the linearized Hamiltonian
system along the brake orbit $\gamma$, normalized by $\Phi(0)=\Id$, and set
$M_\gamma:=\Phi(T)\in\Sp(2n)$. Recall that in the notation of \cite{HWY20} one has
$m^-(x,\Lambda_P)=\iMorT{\gamma}$ for the associated brake orbit $x$ with periodic
boundary condition $\Lambda_P$.

Assume by contradiction that $\gamma$ is linearly stable in the sense of
\cite{HWY20}, i.e.\ $\sigma(M_\gamma)\subset\U$ and $M_\gamma$ is semisimple.
Let $\C_-:=\{\lambda\in\C\mid \Im(\lambda)<0\}$. By \cite[Theorem~1.9]{HWY20},
\begin{equation}\label{eq:HWY-nonreal-bound}
\dim_{\C}\bigoplus_{\lambda\in\sigma(M_\gamma)\cap\C_-}
\ker\!\bigl(M_\gamma-\lambda\Id_{2n}\bigr)
\le 2\,\iMorT{\gamma}.
\end{equation}

Since $M_\gamma$ is real symplectic and $\sigma(M_\gamma)\subset\U$, every
non-real eigenvalue $\lambda\in\U\setminus\R$ comes with its complex conjugate
$\bar\lambda$, and the corresponding eigenspaces have the same complex
dimension. Therefore,
\[
\dim_{\C}\bigoplus_{\lambda\in\sigma(M_\gamma)\setminus\R}\ker(M_\gamma-\lambda\Id_{2n})
=
2\,
\dim_{\C}\bigoplus_{\lambda\in\sigma(M_\gamma)\cap\C_-}\ker(M_\gamma-\lambda\Id_{2n})
\le 4\,\iMorT{\gamma},
\]
where we used \eqref{eq:HWY-nonreal-bound} in the last step. Since $M_\gamma$ is
semisimple, the direct sum of its eigenspaces equals $\C^{2n}$, hence
\begin{equation}\label{eq:real-eigenspace-lowerbound}
\dim_{\C}\bigoplus_{\lambda\in\sigma(M_\gamma)\cap\R}\ker(M_\gamma-\lambda\Id_{2n})
\ge 2n-4\,\iMorT{\gamma}
=2\bigl[n-2\,\iMorT{\gamma}\bigr].
\end{equation}
By the hypothesis $n-2\,\iMorT{\gamma}\ge 1$, the right-hand side is at least $2$.

Under $\sigma(M_\gamma)\subset\U$, the only real points on $\U$ are $\lambda=\pm1$.
Consequently \eqref{eq:real-eigenspace-lowerbound} implies
\[
\dim\ker(M_\gamma-\Id_{2n})+\dim\ker(M_\gamma+\Id_{2n})
\ge 2\bigl[n-2\,\iMorT{\gamma}\bigr]\ge 2,
\]
so at least one of $\pm1$ occurs as a Floquet multiplier.

If $\gamma$ is \emph{strongly nondegenerate}, i.e.\ $\pm1\notin\sigma(M_\gamma)$,
this is a contradiction. Therefore $\gamma$ cannot be linearly stable.
\qed


\subsection*{Proof of Corollary~\ref{thm:main-3-intro}}

\begin{proof}
The first claim is immediate from Theorem~\ref{thm:main-2-intro}. Indeed, under
the present assumptions $\dim M\ge 3$ and $\gamma$ is a non-degenerate
mountain--pass brake solution. In particular $\gamma$ is non-constant and, by
the mountain--pass characterization, the fixed--time Morse index satisfies
$\iMorT{\gamma}\le 1$. On the other hand, by Theorem~\ref{thm:main-1-intro} one
has $\iMorT{\gamma}\ge 1$. Hence
\[
\iMorT{\gamma}=1.
\]
Therefore $\dim M-2\,\iMorT{\gamma}=\dim M-2\ge 1$, and Theorem~\ref{thm:main-2-intro}
yields that $\gamma$ is not linearly stable.

For the second claim, assume that $M_\gamma$ is semisimple. If $\gamma$ were
spectrally stable, then $\sigma(M_\gamma)\subset\U$, and together with
semisimplicity this would imply that $\gamma$ is linearly stable in the sense of
\cite{HWY20}. This contradicts the first part. Hence $\gamma$ cannot be
spectrally stable, i.e.\ it is spectrally unstable.
\end{proof}


\section{Applications}

We compute fixed-time and free-time Morse indices for: (i) the planar (an)isotropic oscillator, (ii) the planar pendulum, and (iii) the planar Kepler problem.

\subsection{The (an)isotropic oscillator}\label{sec:AO}

Consider
\[
L(x,y,\dot x,\dot y)
=\frac12(\dot x^2+\dot y^2)
-\left(\frac{x^2}{2}+\mu^2\frac{y^2}{2}\right),
\]
with Euler--Lagrange equations
\begin{equation}\label{eq:EL-HO}
\ddot x + x =0,\qquad 
\ddot y + \mu^2 y = 0.
\end{equation}
For $k>0$ the Hill region is
\[
\mathcal H_k
=\Bigl\{(x,y)\in\R^2 \mid \tfrac{x^2}{2}+\mu^2\tfrac{y^2}{2}\le k\Bigr\},
\qquad
\partial\mathcal H_k
=\Bigl\{\tfrac{x^2}{2}+\mu^2\tfrac{y^2}{2}=k\Bigr\}.
\]

With $(p_x,p_y)=(\dot x,\dot y)$,
\[
H(p_x,p_y,x,y)
=\frac12(p_x^2+p_y^2)+\frac{x^2}{2}+\mu^2\frac{y^2}{2},
\]
and
\begin{equation}\label{eq:Ham-HO}
\dot p_x=-x,\qquad 
\dot p_y=-\mu^2 y,\qquad
\dot x=p_x,\qquad 
\dot y=p_y.
\end{equation}
Let $\mathcal S(p_x,p_y,x,y)=(-p_x,-p_y,x,y)$; then $\Fix(\mathcal S)=\{p_x=0,\ p_y=0\}$ and a brake orbit satisfies $z(-t)=\mathcal S(z(t))$.

For $k>0$ define the vertical brake orbit
\[
x_k(t)\equiv0,\qquad 
y_k(t)=a_k\cos(\mu t),\qquad 
a_k=\sqrt{\frac{2k}{\mu^2}},
\]
so
\[
z_k(t)=(0,-\mu a_k\sin(\mu t),0,a_k\cos(\mu t))
\]
has period $T_k=\frac{2\pi}{\mu}$ and satisfies $z_k(-t)=\mathcal S(z_k(t))$. The turning points ($p_y=\dot y=0$) occur at $t=0$ and $t=T_k/2$.

\subsection*{Fixed-time index}

The index forms are
\[
Q_x[\xi]
=\int_0^{T_k}\bigl(\xi'^2-\xi^2\bigr)\,dt,
\qquad
Q_y[\eta]
=\int_0^{T_k}\bigl(\eta'^2-\mu^2\eta^2\bigr)\,dt,
\]
with periodic boundary conditions. The Fourier modes give eigenvalues
\[
\lambda_n=n^2\mu^2-1
\quad\text{for }-\frac{d^2}{dt^2}-1,
\qquad
\lambda_n=(n^2-1)\mu^2
\quad\text{for }-\frac{d^2}{dt^2}-\mu^2,
\quad n\in\Z.
\]
Counting real modes yields
\begin{equation}\label{eq:AO-fixed-time-index}
\iMorT{\gamma_k}
=2\Big\lceil\frac{1}{\mu}\Big\rceil.
\end{equation}

\subsection*{Dirichlet index}

For the Dirichlet splitting on $[0,T_k]$,
\[
L_1=-\frac{d^2}{dt^2}-1,
\qquad
L_2=-\frac{d^2}{dt^2}-\mu^2,
\qquad
\eta(0)=\eta(T_k)=0,
\]
we have $\iMorD{\gamma_k}=\iMorD{L_1}+\iMorD{L_2}$. The eigenvalues
\[
\lambda^{(1)}_n=\Bigl(\frac{\pi n}{T_k}\Bigr)^2-1,
\qquad
\lambda^{(2)}_n=\Bigl(\frac{\pi}{T_k}\Bigr)^2(n^2-4),
\quad n\in\N,
\]
give
\begin{equation}\label{eq:Dirichlet-index-L1}
\iMorD{L_1}
=\Bigl\lceil\frac{2}{\mu}\Bigr\rceil-1,
\end{equation}
\begin{equation}\label{eq:Dirichlet-index-L2}
\iMorD{L_2}=1,
\qquad
\dim\ker(L_2)=1,
\end{equation}
and therefore
\begin{equation}\label{eq:Dirichlet-total-two-blocks}
\iMorD{\gamma_k}
=\Bigl\lceil\frac{2}{\mu}\Bigr\rceil
=\Bigl\lceil\frac{T_k}{\pi}\Bigr\rceil.
\end{equation}

\subsection*{Free-time index}

Along the orbit cylinder $\ddot y+\mu^2y=0$ with energy $\frac12(\dot y^{\,2}+\mu^2y^2)=k$, one has $T(k)=\frac{2\pi}{\mu}$, so $T'(k)=0\ge0$. Hence the period direction contributes $+1$ and
\[
\iMorE{\gamma_k}
=2\Big\lceil\frac{1}{\mu}\Big\rceil+1.
\]

\subsection{The planar pendulum: period and indices}
\label{sec:pendulum-morse-separatrix}

In canonical coordinates $(p,q)\in\R^2$ consider
\[
H(p,q)=\frac12p^2+(1-\cos q),
\qquad
J=\begin{pmatrix}0&-1\\ 1&0\end{pmatrix},
\]
so Hamilton's equations are
\begin{equation}\label{eq:Ham-pendulum-merged}
\dot p=-\sin q,\qquad \dot q=p.
\end{equation}
For $k>0$, $\mathcal H_k=\{q\in\R:1-\cos q\le k\}$. In the libration regime $k\in(0,2)$ one has $\mathcal H_k=(-a_k,a_k)$ with $a_k=\arccos(1-k)\in(0,\pi)$ and $\partial\mathcal H_k=\{-a_k,a_k\}$. For each $k\in(0,2)$ there is a unique (up to time shift) librational periodic orbit $z_k(t)=(p_k(t),q_k(t))\subset H^{-1}(k)$, normalized by
\begin{equation}\label{eq:turning-phase}
z_k(0)=(0,a_k),
\end{equation}
with minimal period $T_k$.


The system is reversible under $\mathcal R(p,q)=(-p,q)$. With \eqref{eq:turning-phase},
\[
z_k(-t)=\mathcal R(z_k(t)),
\]
so $z_k$ is a brake orbit with brake instants $t=0$ and $t=T_k/2$.

\subsection*{Monotonicity of the period}

Energy conservation gives
\[
\frac12\dot q(t)^2+(1-\cos q(t))=k,
\qquad
\dot q(t)=\pm\sqrt{2\bigl(k-(1-\cos q)\bigr)},
\]
hence
\[
T_k=4\int_0^{a_k}\frac{dq}{\sqrt{2\bigl(k-(1-\cos q)\bigr)}}.
\]

\begin{lem}[Strict increase of the libration period]\label{lem:Tk-increasing-merged}
For $k\in(0,2)$, $T_k$ is smooth and
\[
\frac{dT_k}{dk}>0,\qquad T_k\to+\infty \ \text{as } k\to2^-.
\]
\end{lem}

\begin{proof}
Using $1-\cos q=2\sin^2(q/2)$ and $q=2\arcsin(\sqrt{k/2}\sin\theta)$,
\[
T_k=4\int_0^{\pi/2}\frac{d\theta}{\sqrt{1-\frac{k}{2}\sin^2\theta}}
=4K\!\left(\frac{k}{2}\right),
\]
where $K$ is the complete elliptic integral of the first kind. Since $K'(m)>0$ for $m\in(0,1)$, $dT_k/dk>0$, and $K(k/2)\to+\infty$ as $k\to2^-$.
\end{proof}


\subsection*{Variational equation, monodromy, and Jordan form}

Let $\Phi_k(t)\in\Sp(2)$ solve
\begin{equation}\label{eq:var-eq-merged}
\dot\Phi_k=J D^2H(z_k(t))\,\Phi_k,
\qquad
\Phi_k(0)=I,
\end{equation}
and set $M_k=\Phi_k(T_k)$. Since $D^2H(p,q)=\mathrm{diag}(1,\cos q)$, the linearized equation is
\[
\dot{\delta p}=-\cos(q_k(t))\,\delta q,\qquad
\dot{\delta q}=\delta p,
\]
equivalently $\ddot{\delta q}+\cos(q_k(t))\delta q=0$.

The vectors $u_k(t)=\dot z_k(t)$ and $v_k(t)=\partial_k z_k(t)$ solve \eqref{eq:var-eq-merged}. Differentiating $z_k(t+T_k)=z_k(t)$ with respect to $k$ gives
\begin{equation}\label{eq:v-monodromy-merged}
v_k(T_k)=v_k(0)-T_k' u_k(0),
\end{equation}
hence
\[
M_k u_k(0)=u_k(0),
\qquad
M_k v_k(0)=v_k(0)-T_k' u_k(0).
\]

At $t=0$, $u_k(0)=(-\sin a_k,0)$ and $v_k(0)=(0,a_k')$ with $a_k'=1/\sin a_k>0$. Since the symplectic form satisfies $\omega(u_k(0),v_k(0))=-1$, the basis $(u_k(0),-v_k(0))$ is symplectic and
\begin{equation}\label{eq:Jordan-positive-merged}
[M_k]_{(u_k(0),-v_k(0))}
=
\begin{pmatrix}1&T_k'\\ 0&1\end{pmatrix}
=:N_1(1,b_k),
\qquad b_k=T_k'>0.
\end{equation}


\subsection*{Maslov/CLM index from the shear normal form}

Let $\Delta=\Graph\Id$. Since $\nu_1(M_k)=1$ for all $k\in(0,2)$, no eigenvalue
crosses $1$ along the family, and the CLM index
\[
\iCLM(\Delta,\Phi_k(t),t\in[0,T_k])
\]
is locally constant in $k$.

For $k>0$ sufficiently small, $z_k$ converges in $C^1$ on compact intervals to
the solution $z_0$ of the linear oscillator
\[
\dot p=-q,\qquad \dot q=p,
\]
whose fundamental solution is $\Psi(t)=e^{Jt}$ and satisfies $\Psi(2\pi)=\Id$.
Moreover $T_k\to2\pi$ as $k\downarrow0$.
\begin{lem}
    Under the above notation, the following holds: 
\begin{equation}\label{eq:Phi-to-Psi-on-fixed}
\sup_{t\in[0,2\pi]}\|\Phi_k(t)-\Psi(t)\|\xrightarrow[k\downarrow 0]{}0.
\end{equation}
\end{lem}
\begin{proof}
    The proof directly follows by continuous dependence of the fundamental solutions, since $D^2H(z_k(t)) \to \Id$ in $\mathscr C^0$. 
\end{proof}
\begin{lem}\label{lem:Phi-tildePsi-C0}
Set $M_k:=\Phi_k(T_k)$ and define $\widetilde\Psi\colon[0,T_k]\to\R^{2n\times 2n}$ by
\[
\widetilde\Psi(t)=
\begin{cases}
\Psi(t), & t\in[0,2\pi],\\[6pt]
\displaystyle
\frac{T_k-t}{T_k-2\pi}\,\Psi(2\pi)+\frac{t-2\pi}{T_k-2\pi}\,M_k,
& t\in[2\pi,T_k].
\end{cases}
\]
Then
\[
\|\Phi_k-\widetilde\Psi\|_{\mathscr C^0([0,T_k])}\xrightarrow[k\downarrow 0]{}0.
\]
\end{lem}

\begin{proof}
On $[0,2\pi]$ we have $\widetilde\Psi=\Psi$, hence by Equation~\eqref{eq:Phi-to-Psi-on-fixed}, we get
\begin{equation}\label{eq:C0-first-piece}
\sup_{t\in[0,2\pi]}\|\Phi_k(t)-\widetilde\Psi(t)\|\to 0.
\end{equation}

It remains to control the difference on the short interval $[2\pi,T_k]$ (which
shrinks as $k\downarrow 0$). For $t\in[2\pi,T_k]$ we write
\[
\Phi_k(t)-\widetilde\Psi(t)
=\bigl(\Phi_k(t)-\Phi_k(2\pi)\bigr)
+\bigl(\Phi_k(2\pi)-\Psi(2\pi)\bigr)
+\bigl(\Psi(2\pi)-\widetilde\Psi(t)\bigr).
\]
Taking norms and the supremum over $t\in[2\pi,T_k]$ gives
\begin{equation}\label{eq:C0-tail-split}
\sup_{t\in[2\pi,T_k]}\|\Phi_k(t)-\widetilde\Psi(t)\|
\le A_k+B_k+C_k,
\end{equation}
where
\[
A_k:=\sup_{t\in[2\pi,T_k]}\|\Phi_k(t)-\Phi_k(2\pi)\|,
\quad
B_k:=\|\Phi_k(2\pi)-\Psi(2\pi)\|,
\quad
C_k:=\sup_{t\in[2\pi,T_k]}\|\Psi(2\pi)-\widetilde\Psi(t)\|.
\]
By \eqref{eq:Phi-to-Psi-on-fixed}, we have $B_k\to 0$ and in particular
$\Phi_k(2\pi)\to\Psi(2\pi)$.

Next, $\widetilde\Psi$ interpolates linearly between $\Psi(2\pi)$ and $M_k$, so
\[
C_k\le \|M_k-\Psi(2\pi)\|=\|\Phi_k(T_k)-\Psi(2\pi)\|.
\]
Since $T_k\to 2\pi$ and each $\Phi_k$ is continuous, we have
$\|\Phi_k(T_k)-\Phi_k(2\pi)\|\to 0$ for each fixed $k$. Moreover, the family
$\{\Phi_k\}$ is uniformly continuous near $t=2\pi$ for $k$ small (indeed,
\eqref{eq:Phi-to-Psi-on-fixed} implies uniform closeness to the continuous map
$\Psi$, hence equicontinuity on $[0,2\pi]$). Therefore $A_k\to 0$.

Finally,
\[
\|M_k-\Psi(2\pi)\|
\le \|\Phi_k(T_k)-\Phi_k(2\pi)\|+\|\Phi_k(2\pi)-\Psi(2\pi)\|
\to 0,
\]
so also $C_k\to 0$. Plugging $A_k,B_k,C_k\to 0$ into \eqref{eq:C0-tail-split}
yields
\[
\sup_{t\in[2\pi,T_k]}\|\Phi_k(t)-\widetilde\Psi(t)\|\to 0.
\]
Together with \eqref{eq:C0-first-piece} this proves the claim.
\end{proof}
By homotopy invariance of the CLM index with fixed endpoints,
\[
\iCLM(\Delta,\Phi_k(t),t\in[0,T_k])
=
\iCLM(\Delta,\widetilde\Psi(t),t\in[0,T_k]).
\]
Since $\widetilde\Psi$ has no crossings with $\Delta$ and even with $\Lambda_D$ on $[2\pi,T_k]$,
\begin{multline}
\iCLM(\Delta,\Phi_k(t),t\in[0,T_k])
=
\iCLM(\Delta,\Psi(t),t\in[0,2\pi])
=2\\[5pt]
\iCLM(\Lambda_D,\Phi_k(t),t\in[0,T_k])
=
\iCLM(\Lambda_D,\Psi(t),t\in[0,2\pi])
=2
\end{multline}
By the Morse–Maslov correspondence for brake-symmetric variations, this yields
\[
\iMorT{z_k}=1.
\]

\begin{thm}[Fixed-time and free-time  Morse indices of pendulum brake orbits]
\label{thm:pendulum-fixed-time-energy}
Let $z_k$ be the librational brake orbit of energy $k\in(0,2)$ and minimal period
$T_k$.  Then:
\[
\iMorT{z_k}=1 \qquad\textrm { and } \qquad 
\iMorE{z_k}=2, \qquad\textrm { and } \qquad \iMorD{z_k}=1. 
\]
In particular,  librational brake orbits  are
minimizers for the fixed-time or free-time  variational problems.
\end{thm}

\begin{proof}
The first equality is proved above (librations by the shear/CLM computation and stability of the crossing count on compact
intervals). For the second, note that $T_k'>0$ by Lemma~\ref{lem:Tk-increasing-merged},
so the orbit-cylinder correction gives $\iMorE{z_k}=\iMorT{z_k}+1=2$ for all
$k\in(0,2)$. For the separatrix, one can either pass to the limit in the
free-period index (using the same compact-interval approximation and the fact
that the extra $+1$ comes from the period direction) or appeal to a direct
argument as in \cite[Proposition~5.1]{HPWX25}. In either case,
$\iMorE{z_{\mathrm{sep}}}=2$.
\end{proof}

\begin{rmk}
For a direct proof of the separatrix (heteroclinic) fixed-time index, see \cite[Proposition~5.1]{HPWX25}.
\end{rmk}
\subsection{The planar Kepler problem}\label{sec:kepler-morse}

On $Q=\R^2\setminus\{0\}$ consider $H(p,q)=\frac12|p|^2-\frac1{|q|}$. For energy $E$, $\mathcal H_E=\{q\in Q: E+\frac1{|q|}\ge0\}$. If $E<0$, then $\mathcal H_E=\{r\le r_E\}$ with $r_E=(-E)^{-1}$ and $\partial\mathcal H_E=\{r=r_E\}$; we restrict to $E<0$.

In polar coordinates $(r,\theta)$ with momenta $(p_r,\ell)$,
\begin{equation}\label{eq:kepler-polar}
H=\frac12\Bigl(p_r^2+\frac{\ell^2}{r^2}\Bigr)-\frac1r,\qquad
\dot r=p_r,\qquad \dot\theta=\frac{\ell}{r^2}.
\end{equation}
The system is reversible under $\mathcal R(p,q)=(-p,q)$; a brake orbit satisfies $z(-t)=\mathcal R(z(t))$.

\subsection*{Elliptic orbits}

For $E<0$ and $\ell\neq0$ the motion is elliptic. Fixing the phase at apocenter ($\dot r(0)=0$) yields a $T(E)$--periodic brake orbit. Kepler’s law gives
\begin{equation}\label{eq:kepler-period}
T(E)=2\pi(-2E)^{-3/2},\qquad E<0,
\end{equation}
hence $\frac{dT}{dE}(E)>0$. Collision-free ellipses are strong local minimizers of the fixed-time action in their homotopy class, so
\begin{equation}\label{eq:kepler-fixed-time-zero}
\iMorT{z_E}=0.
\end{equation}
Since $E\mapsto(z_E,T(E))$ is a smooth orbit cylinder with $T'(E)>0$,
\begin{equation}\label{eq:kepler-free-time}
\iMorE{z_E}=\iMorT{z_E}+1=1.
\end{equation}

\subsection*{The ejection--collision brake orbit}

For $\ell=0$ the motion is radial and satisfies
\begin{equation}\label{eq:radial-kepler-final}
E=\frac12\dot r^{\,2}-\frac1r,\qquad \ddot r=-\frac1{r^2}.
\end{equation}
With $r(0)=r_E=(-E)^{-1}$ and $\dot r(0)=0$, the solution reaches collision in finite time; by time reversal we obtain the ejection--collision brake orbit $z_{\mathrm{ec}}$. The flight time is
\begin{equation}\label{eq:kepler-flight-final}
t_c(E)=\int_0^{r_E}\frac{dr}{\sqrt{2(E+\frac1r)}}=\frac{\pi}{2}(-2E)^{-3/2},
\qquad
T(E)=2t_c(E)=\pi(-2E)^{-3/2},
\end{equation}
so $T'(E)>0$.

\subsubsection*{Levi--Civita regularization}

Identify $\R^2\simeq\C$ and set $q=u^2$, $r=|u|^2$. The cotangent lift defined by $p\,dq=v\,du$ gives
\begin{equation}\label{eq:LC-cotangent-lift-final}
v=2u\,p,\qquad p=\frac{v}{2u},
\end{equation}
so $\LC(u,v)=(q,p)$ is symplectic for $u\neq0$. With time change $dt/ds=|u|^2$, the Poincar\'e Hamiltonian is
\begin{equation}\label{eq:Poincare-K-final}
K(u,v)=\frac18|v|^2-E|u|^2-1.
\end{equation}
Along the radial orbit ($u,v\in\R$) the equation reduces to
\begin{equation}\label{eq:oscillator-final}
u''+\omega_E^2u=0,\qquad \omega_E=\sqrt{-2E},
\end{equation}
and since $q=u^2$, the physical arc closes after half a period,
\begin{equation}\label{eq:regularized-period-final}
S(E)=\frac{\pi}{\omega_E}.
\end{equation}

\subsubsection*{Fixed- and free-time indices}

On $[0,S(E)]$ the second variation (up to a positive factor) is
\begin{equation}\label{eq:Q-osc-final}
\mathsf Q_E(\eta)=\int_0^{S(E)}\bigl(\eta'^2-\omega_E^2\eta^2\bigr)\,ds,
\end{equation}
with periodic boundary conditions. For $\eta_m(s)=e^{i\frac{2\pi m}{S(E)}s}$,
\[
\lambda_m=\Bigl(\frac{2\pi m}{S(E)}\Bigr)^2-\omega_E^2=(4m^2-1)\omega_E^2,
\]
so only $\lambda_0<0$ and
\begin{equation}\label{eq:kepler-ec-fixed-time-final}
\iMorT{z_{\mathrm{ec}}}=1.
\end{equation}
Since $T'(E)>0$,
\begin{equation}\label{eq:kepler-ec-free-time-final}
\iMorE{z_{\mathrm{ec}}}=\iMorT{z_{\mathrm{ec}}}+1=2.
\end{equation}

\begin{thm}
For $E<0$, the radial ejection--collision brake trajectory satisfies
$\iMorT{z_{\mathrm{ec}}}=1$ and $\iMorE{z_{\mathrm{ec}}}=2$.
\end{thm}

\begin{rmk}
The Levi--Civita map \eqref{eq:LC-cotangent-lift-final} is symplectic on $u\neq0$, the time change is encoded in \eqref{eq:Poincare-K-final}, and collision becomes regular; the index reduces to the Fourier computation of \eqref{eq:oscillator-final}.
\end{rmk}

\begin{rmk}[Keplerian ellipses]
For $\ell\neq0$, elliptic solutions satisfy $\iMorT{z_E}=0$ and, since $T'(E)>0$, $\iMorE{z_E}=1$.
\end{rmk}



\vskip0.7cm

\noindent
\textsc{JProf. Dr. Luca Asselle}\\
Ruhr-Universit\"at Bochum\\
Fakult\"at f\"ur Mathematik\\
Universit\"atsstra{\ss}e 150\\
44801 Bochum, Germany\\
\texttt{luca.asselle@rub.de}

\vskip0.6cm

\noindent
\textsc{Prof. Dr. Xijun Hu}\\
School of Mathematics\\
Shandong University\\
Jinan 250100, P. R. China\\
\texttt{xjhu@sdu.edu.cn}

\vskip0.6cm

\noindent
\textsc{Prof. Dr. Alessandro Portaluri}\\
Universit\`a degli Studi di Torino (DISAFA)\\
Largo Paolo Braccini 2\\
10095 Grugliasco (Torino), Italy\\
\texttt{alessandro.portaluri@unito.it}

\medskip
\noindent
Visiting Professor of Mathematics\\
New York University Abu Dhabi\\
Saadiyat Marina District, Abu Dhabi, UAE\\
\texttt{ap9453@nyu.edu}

\vskip0.6cm

\noindent
\textsc{Prof. Dr. Li Wu}\\
School of Mathematics\\
Shandong University\\
Jinan 250100, P. R. China\\
\texttt{vvvli@sdu.edu.cn}

\end{document}